 \documentclass[12pt]{article}
 \usepackage{amsmath,amsfonts,amssymb,theorem}
 \numberwithin{equation}{section}

 \newtheorem{prop}{Proposition}[section]

 \newtheorem{cor}{Corollary}[section]

 \newtheorem{thm}{Theorem}[section]
 \newtheorem{lem}{Lemma}[section]

 \theorembodyfont{\rmfamily}
 \newtheorem{dfn}{Definition}[section]
 \newtheorem{mcdn}{Main Condition}[section]
 \newtheorem{warning}{Warning}
 
 \newtheorem{rmk}{Remark}

 \newenvironment{pf}{\paragraph{Proof}}{\par\bigskip}
 \newcommand{\qed}{\ifhmode\unskip\nobreak\fi\quad\ensuremath\square}

\newenvironment{digr}[1]{\subsubsection*{#1}}{\quad$\diamondsuit$\par\bigskip}

 \newcommand{\Span}[1]{\left< #1 \right>}
 \newcommand{\half}{\textstyle\frac12}
 \newcommand{\dd}{\mathrm{d}}
 \newcommand{\id}{\operatorname{id}}
 \newcommand{\pr}{\operatorname{pr}}
 \newcommand{\ch}{\operatorname{ch}}
 \newcommand{\td}{\operatorname{td}}
 \newcommand{\CD}{\operatorname{CD}}
 \newcommand{\KH}{\mathrm{KH}}
 \newcommand{\LC}{\mathrm{LC}}
 \newcommand{\SDAG}{\mathrm{SDAG}}
 \newcommand{\II}{\mathrm{II}}
 \newcommand{\gen}{\mathrm{gen}}
 \newcommand{\rest}[1]{_{{\textstyle{|}}#1}} 

 \newcommand{\dbar}{\overline{\partial}}
 \newcommand{\ubar}{\overline{u}}
 \newcommand{\Ibar}{\overline{I}}
 \newcommand{\iso}{\cong}
 \newcommand{\wave}{\widetilde}
 \newcommand{\tensor}{\otimes}
 \newcommand{\ov}{\overline}

 \newcommand{\sC}{\mathcal C}
 \newcommand{\sF}{\mathcal F}
 \newcommand{\sH}{\mathcal H}
 \newcommand{\sK}{\mathcal K}
 \newcommand{\sL}{\mathcal L}
 \newcommand{\sM}{\mathcal M}
 \newcommand{\SM}{\mathcal{SM}}
 \newcommand{\Oh}{\mathcal O}
 \newcommand{\sR}{\mathcal R}

 \newcommand{\fie}{\varphi}
 \newcommand{\ga}{\gamma}
 \newcommand{\ka}{\kappa}
 \newcommand{\om}{\omega}

 \newcommand{\Ga}{\Gamma}
 \newcommand{\La}{\Lambda}
 \newcommand{\Om}{\Omega}
 \newcommand{\la}{\lambda}
 \newcommand{\Si}{\Sigma} 

 \newcommand{\PP}{\mathbb P}
 \newcommand{\C}{\mathbb C}
 \newcommand{\HH}{\mathbb H}
 \newcommand{\Q}{\mathbb Q}
 \newcommand{\R}{\mathbb R}
 \newcommand{\Z}{\mathbb Z}

 \newcommand{\ad}{\operatorname{ad}} 
 \newcommand{\rk}{\operatorname{rank}}
 \newcommand{\hcf}{\operatorname{hcf}}
 \newcommand{\mir}{\operatorname{mir}}
 \newcommand{\topo}{\operatorname{top}}
 \newcommand{\vdim}{\operatorname{v.dim}}
 \newcommand{\Pic}{\operatorname{Pic}} 
 \newcommand{\End}{\operatorname{End}}
 \newcommand{\GFT}{\operatorname{GFT}}
 \newcommand{\sGFT}{\operatorname{sGFT}}
 \newcommand{\LaGr}{\operatorname{\La_{\uparrow}\!}}
 \newcommand{\Hom}{\operatorname{Hom}}
 \renewcommand{\Im}{\operatorname{Im}}
 \renewcommand{\Re}{\operatorname{Re}}
 \newcommand{\Sing}{\operatorname{Sing}}
 \newcommand{\Vol}{\operatorname{Vol}}

 \newcommand{\SO}{\operatorname{SO}}
 \newcommand{\U}{\operatorname U} 
 \newcommand{\pt}{\mathrm{pt}}

\begin{document}

 \title{Geometric quantization and mirror symmetry}
 \markright{\qquad Geometric quantization and mirror symmetry}

 \author{Andrei Tyurin}
 \date{}
 \maketitle
 \begin{abstract}
After the appearance of \cite{T3} I received an e-mail from Cumrun Vafa,
who recognized that the subject is closely related to that of his preprint
\cite{V}. This text started out as an e-mail ``reply'' to his letter. All
the constructions we propose have well known ``spectral curve'' prototypes
(see for example Friedman and other \cite{FMW}, Bershadsky and other
\cite{BJPS} and a number of others). Roughly speaking, our constructions
are the {\em spectral curve} construction plus the {\em phase geometry}
described in \cite{T3}. So this text should really come before \cite{T3},
as motivation for the development of the geometry of the phase map in
\cite{T3}. 
 \end{abstract}

 \section{spLag cycles}
 \markright{\qquad Geometric quantization and mirror symmetry}

We begin by recalling the actual geometric construction for a pair
$\sL\subset S$, where $S$ is a smooth symplectic manifold of dimension
$2n$ with a given tame almost complex structure $I$, and $\sL\subset S$ a
smooth, oriented Lagrangian submanifold (of maximal dimension
$\dim\sL=n=\half\dim S$); this construction has recently become quite
popular in the set-up of Calabi--Yau threefolds. The structure on $S$ is
an {\em almost K\"ahler structure}, and we say for short that $S$ is an
{\em aK manifold}. Write $\om$ for the symplectic form on $S$ and $I$ for
the almost complex structure, so that the tangent space $TS_p$ at a point
$p$ is $\C^n$ with the constant symplectic form $\Span{\,\,,\,}=\om_p$ and
the constant Euclidean metric $g_p$, giving the Hermitian triple
$(\om_p,I_p,g_p)$.

We now define the {\em Lagrangian Grassmannian} $\LaGr p=\LaGr(TS_p)$ to
be the Grassmannian of maximal oriented Lagrangian subspaces in $TS_p$.
Taking this space over every point of $S$ gives the {\em oriented
Lagrangian Grassmannization} of $TS$:
 \begin{equation}
 \pi\colon \LaGr(S)\to S \quad\text{with} \quad\pi^{-1} (p) =\LaGr p.
 \label{eq1.1}
 \end{equation}

Our tame almost complex structure on $S$ gives each fibre the standard
form
 \begin{equation}
 \LaGr p=\U(n)/\SO(n)
 \label{eq1.2}
 \end{equation}
This space admits a canonical map
 \begin{equation}
 \det\colon\LaGr p\to\U(1)=S^1_p
 \quad\text{sending $u\in\U(n)$ to $\det u\in\U(1)=S^1$.}
 \label{eq1.3}
 \end{equation}
Recall that the inverse image of the fundamental class of $S^1$ on $\LaGr
p$ is the {\em universal Maslov class}. Taking this map over every point
of $S$ gives a map 
 \begin{equation} 
\det\colon \LaGr(S)\to S^1(L_{-K}),
 \label{eq1.4}
 \end{equation}
where $S^1(L_{-K})$ is the unit circle bundle of the line bundle
$\bigwedge^n TS=\det TS$, with first Chern class
 \begin{equation}
 c_1(\det TS)=-K_S,
 \label{eq1.5}
 \end{equation}
where $K_S$ is the canonical class of $S$. Recall that, as a cohomology
class, $K_S$ does not depend on the compatible almost complex structure.

Now for every oriented Lagrangian cycle $\sL\subset S$, we have
the Gauss lift of the embedding $i\colon\sL\to S$ to a section 
 \begin{equation}
 G(i)\colon\sL\to\LaGr(S)\rest{\sL},
 \label{eq1.6}
 \end{equation} sending a point $p\in\sL$
to the oriented subspace $T\sL_p\subset TS_p$. The composite of this
Gauss map with the projection (\ref{eq1.4}) gives a map 
 \begin{equation}
 \det \circ G(i)\colon\sL\to S^1(L_{-K})\rest{\sL}.
 \label{eq1.7}
 \end{equation}

Now suppose that the cohomology class of the symplectic form is
proportional to the canonical class of $S$, that is,
 \begin{equation}\ka\cdot[\om]=K_S \quad\text{for some $\ka\in \Q$};
 \label{eq1.8}
 \end{equation}
then the restriction $\det TS\rest{\sL}$ is topologically trivial,
because the restriction of $[\om]$ to a Lagrangian $\sL$ is zero.

Now for $[\om]\in H^2(S,\Z)$, there exists a complex line bundle $L$ such
that $c_1(L)=[\om]$. Moreover, suppose that $L$ has a Hermitian connection
$a$ with curvature form $F_a=\om$. A quadruple 
 \begin{equation}
 (S,\om,L,a)
 \label{eq1.9}
 \end{equation} 
of this form is called a {\em prequantization} of the classical mechanical
system with phase space $(S,\om)$.

But in our situation (\ref{eq1.8}), the canonical line bundle $K_S$ has a
Hermitian connection $a_K$ with curvature form
 \begin{equation}
 \frac{1}{2\pi i}F_{a_K}=\ka\om.
 \label{eq1.10}
 \end{equation} 
Recall that any curvature form has its coefficients in the Lie algebra of
the gauge group. For our $\U(1)$ case, this Lie algebra is $2\pi i\R$, so
to get a real (integral) form we must multiply the curvature form by
$\frac{1}{2\pi i}$, as in (\ref{eq1.10}). 

Both of these connections restrict to flat connections on $\sL$, or
equivalently, to characters
 \[
 \chi\colon\pi_1(\sL)\to\U(1)
 \]
and $\chi^{\ka}$ of the fundamental group.

Topologically, there exists a trivialization 
 \[
 S^1(L_{-K})\rest{\sL}=\sL\times S^1,
 \]
but we would like to construct a canonical projection
 \begin{equation}
 \pr\colon S^1(L_{-K})\rest{\sL}\to S^1
 \label{eq1.11}
 \end{equation}
preserving the Hermitian form. Then composing the maps (\ref{eq1.4}),
(\ref{eq1.6}) and (\ref{eq1.11}) would give a map 
 \begin{equation}
 m=\pr\circ\det\circ G(i)\colon\sL\to S^1
 \label{eq1.12}
 \end{equation}
But a priori we cannot do this if the character $\chi^{\ka}$ is not
trivial. To avoid this, we lift these connection to the universal cover
$U(\sL)$ to get the trivial connection with the covariant constant section
 \begin{equation}
 m_I\colon U(\sL)\to\U(1),
 \label{eq1.13}
 \end{equation}
which is a lifting of $m$ (\ref{eq1.12}). We get the functional equation
 \begin{equation}
 m_I(g(u))=\chi^{\ka}(g)\cdot m_I(u).
 \label{eq1.14}
 \end{equation}

 \begin{dfn} This map $m_I$ is called the {\em phase map} with respect to
the almost complex structure $I$.
 \end{dfn}

Thus a priori, the phase map is a multivalued function on a Lagrangian
cycle $\sL$, but its log derivative
 \[
m_I^{-1}\cdot\dd m_I\in \Om^1(\sL)
 \]
is an ordinary differential form on $\sL$.

 \begin{dfn}
 \begin{enumerate}
 \item A cycle $\sL$ is a {\em Bohr--Sommerfeld} cycle of $(S,\om,L,a)$ if
$\chi=1$.

 \item A cycle $\sL$ is called a {\em special Lagrangian cycle} of $S$
({\em spLag cycle} for short) if 
 \[
 m_I^{-1}\cdot\dd m_I=0
 \]
 \end{enumerate}
 \end{dfn}

 \begin{rmk} In this definition, we call $\sL$ a {\em cycle} rather than a
submanifold, because it may be singular. We really only need the Gauss map
(\ref{eq1.6}) to be well defined; thus $\sL$ can have nodes, and so on. Thus
below we call a {\em cycle} any subvariety with a regular Gauss map.
 \end{rmk}

 \begin{digr}{Mirror digression I} (\ref{eq1.8}) and (\ref{eq1.9}) hold
automatically if $S$ is a Calabi--Yau $n$-fold (CY$_n$ for short), that
is, $K_S=0$. Then $\ka=0$ and the canonical line bundle is trivial, with
the trivial connection. In this case, the notion of spLag cycle coincides
with that in calibrated geometry (see Harvey and Lawson \cite{HL}).

Recall that a {\em complex orientation} of a Calabi--Yau manifold $X$ is a
choice of trivialization of the canonical line bundle $L_K$, that is, a
choice of a holomorphic $n$-form $\theta$. For an oriented Calabi--Yau
threefold $(X,\theta)$, a spLag cycle is a 3-dimensional Lagrangian
submanifold $\sL$ such that the restriction of $\theta$ satisfies the
equivalent conditions
 \begin{equation}
 \Im\theta\rest{\sL}=0 \quad \text{and}
 \quad\Re\theta\rest{\sL}=\Vol_g(\sL).
 \label{eq1.15}
 \end{equation}

The local deformation theory of such submanifolds $\sL\subset S$ is well
under\-stood. The tangent space to the moduli space $\sM_{\sL}$ of
deformations at a submanifold $\sL$ is $H^1(\sL,\R)$, as the space of
harmonic 1-forms on $\sL$. This space does not depend on the second
quadratic form or other attributes of embeddings. In particular if
$H^1(\sL,\R)=0$ then $\sL$ is rigid as a special Lagrangian submanifold.
So we can expect that there exists a finite set of such submanifolds
$\{\sL_1,\dots,\sL_N\}$ in one cohomology class $[\sL_i]$. The problem of
computing numbers like $N$ has become quite popular recently, and we would
like to remark that our construction also works for Fano varieties and
varieties of general type.

On the other hand, on a Calabi--Yau manifold, the symplectic form of any
K\"ahler metric $\om$ does not depend on the trivialization of the
canonical class. The cohomology class $[\om]$ is integral and the line
bundle $L$ with $c_1(L)=[\om]$ is holomorphic and carries a Hermitian
connection $a$ with curvature form $\frac{1}{2\pi i}F_a=\om$. So we can
add the Bohr--Sommerfeld condition to the special Lagrangian condition.
Then we expect to get a finite set of spLag Bohr--Sommerfeld cycles in one
cohomology class.
 \end{digr}

 \begin{digr}{Digression on geometric quantization} We recall two
approaches to the quantization of $(S,\om,L,a)$ (\ref{eq1.9}) (see
\'Sniaty\-cki \cite{S} or Woodhouse \cite{W}).

The first method is complex quantization: we give $S$ a complex structure
$I$ such that $S_I$ is a K\"ahler manifold with K\"ahler form $2\pi i\om$.
Then for any {\em level} $k\in\Z^+$, the line bundle $L^k$ is a holomorphic
line bundle on $S_I$. Complex quantization provides the space of {\em wave
functions of level $k$}:
 \begin{equation}
 \sH_{L^k}=H^0(S_I, L^k)
 \label{eq1.16}
 \end{equation}
-- the space of {\em holomorphic} sections of $L^k$. It is reasonable to
call the complex structure $S_I$ and the holomorphic line bundle $L$ a
{\em complex polarization} of $(S,\om,L,a)$ (\ref{eq1.9}). 

Now a real polarization of $(S,\om,L,a)$ is a Lagrangian fibration
 \begin{equation}
 \pi\colon S\to B
 \label{eq1.17}
 \end{equation}
such that every point of $b\in B$
 \[
 \om\rest{\pi^{-1}(b)}=0,
 \]
and for generic $b$, the fibre $\pi^{-1}(b)$ is a smooth Lagrangian. 

Then restricting $L$ to a Lagrangian fibre gives the flat connection
$a\rest{\text{fibre}}$. According to the general theory of real
quantizations, we expect to get a finite number of Bohr--Sommerfeld
Lagrangian cycles among the fibres, and we can construct a new collection
of spaces of wave functions
 \begin{equation}
 \sH_{\pi}^k=\bigoplus_{\text{$k$-BS}} \C\cdot s_{\sL_i},
 \label{eq1.18}
 \end{equation}
where $s_{\sL_i}$ is the covariant constant section of the restriction of
$(L^k, ka)$ to a $k$-Bohr--Sommerfeld fibre $\sL_i$ of the real
polarization $\pi$.

 \begin{rmk} The important observation, proved mathematically in a number
of cases, is that the projectivization of the spaces (\ref{eq1.16}) do not
depend on the choice of complex structure on $S$. They are given purely by
the symplectic prequantization data. The same is true for the
projectivization of the spaces (\ref{eq1.18}). Moreover, these spaces do
not depend on the real polarization $\pi$ (\ref{eq1.18}) (provided that we
extend our prequantization data $(S,\om,L,a,\sF)$ by adding some ``half
density'' $\sF$) (see Guillemin and Sternberg \cite{GS}).
 \end{rmk}

The problem of comparing the ranks
 \[
 \rk \sH_{L^k} \quad \text{and} \quad \rk \sH_{[\sL]}^k
 \]
is the {\em numerical quantization problem}.

In particular, any CY manifold $X$ with any K\"ahler metric can be viewed
as the phase space of a classical mechanical system $(X,\frac{1}{2\pi
i}\om)$, where $\om$ is the K\"ahler form of the metric, with complex
polarization $(L,a)$, such that $F_a=\frac{1}{2\pi i}\om$, as a
prequantization data. Thus complex quantization provides the space of wave
functions (\ref{eq1.16}).

On the other hand, it is easy to see that a smooth spLag Bohr--Sommerfeld
cycle does not have positive dimensional deformations, that is, it is rigid
in its cohomology class $[\sL]$. Then there exists new spaces of wave
functions of any level $k$
 \begin{equation}
 \sH_{[\sL]}^k=\bigoplus_{\text{spLag BS}} \C\cdot s_{\sL_i}
 \label{eq1.19}
 \end{equation}
where $\sL_i$ is a spLag Bohr--Sommerfeld cycle of $(X,\frac{k}{2\pi i}\om,
L^k, ka)$.
 \end{digr}

We now return to the general phase map. The differential of the phase map
has a standard description: let
 \[
\nabla_\LC\colon \Ga (TS)\to \Ga(TS\tensor T^*S)
 \]
be the Levi-Civita connection, which we restrict to the restriction of the
tangent bundle to a Lagrangian cycle $\sL$. Let $N_{\sL\subset S}$ be the
normal bundle of $\sL$ in $S$. Then the Levi-Civita connection defines
connections 
 \[
 \aligned
 \nabla_\LC\colon \Ga (T\sL) &\to \Ga(T\sL\tensor T^*\sL)\\[4pt]
 \nabla_\LC\colon \Ga (N_{\sL\subset S}) &\to \Ga(N_{\sL\subset
S}\tensor T^*\sL)
 \endaligned
 \]
and a tensor
 \[
 \II\colon T\sL\to \Hom(T\sL, N_{\sL\subset S})
 \]
(called the second quadratic form), that is,
 \begin{equation}
\II\in\End T\sL\tensor N_{\sL\subset S}.
 \label{eq1.20}
 \end{equation}
The trace component of $\End T\sL =\C\oplus\ad T\sL$ gives a section of
the normal bundle
 \begin{equation}
H\in \Ga (N_{\sL\subset S}),
 \label{eq1.21}
 \end{equation}
called the {\em mean curvature}.

Let $V$ be a vector field on $\sL$. Then pointwise, the value of the
differential $\dd m_I$ on $V$ is given by the inner product
 \begin{equation}
 \dd m_I(V_p)=(H_p, I(V_p)),
 \label{eq1.22}
 \end{equation}
where $I$ is the operator of our almost complex structure. {From} this we
get immediately:

 \begin{prop} A Lagrangian cycle $\sL$ is spLag iff the mean curvature
vector vanishes at every point $p\in\sL$, that is, $H_p=0$.
 \end{prop}

Thus we see that a spLag cycle is indeed {\em minimal}.

We return to conditions (\ref{eq1.8}--\ref{eq1.9}) in the general situation
(not only CY). For every Lagrangian submanifold $\sL$, restriction defines
a trivial line bundle $(L,a)\rest{\sL}$ with flat connection. In addition
to this, we have the differential 1-form $m_I^{-1}\cdot\dd m_I$.

 \begin{lem}
 \begin{enumerate}
 \item $m_I^{-1}\cdot\dd m_I$ is a closed $1$-form;

 \item its cohomology class is the universal Maslov class of $\sL$.
 \end{enumerate}
 \end{lem}

 \begin{pf} Indeed, this form is the pullback of the constant form on
$S^1=\U(1)$, which is closed and represents the generator of $H^1(S^1,\Z)$.
 \qed \end{pf}

We correct the flat connection $a\rest{\sL}$ by this differential form,
setting:
 \begin{equation}
 a_c=a\rest{\sL}-\frac{1}{\ka}\cdot m_I^{-1}\dd m_I.
 \label{eq1.23}
 \end{equation}
The flat connection $a_c$ is equivalent to a character of the fundamental
group
 \begin{equation}
\chi_c\colon\pi_1(\sL)\to\U(1)
 \label{eq1.24}
 \end{equation}
and we can correct the notion of a Bohr--Sommerfeld cycle:

 \begin{dfn}\label{dfnBS} A Lagrangian cycle $\sL$ is an {\em $m$-corrected
Bohr--Sommerfeld cycle} (cBS cycle) if $\chi_c=1$.
 \end{dfn}

In the general case it is reasonable to expect that for any real
polarization $\pi$ (\ref{eq1.17}), there exists only a finite set of cBS
fibres, and we get the corrected version of the spaces (\ref{eq1.18}):
 \begin{equation}
 \sH_{[\sL]}^k=\bigoplus_{\text{spLag $k$-cBS}} \C\cdot s_{\sL_i},
 \label{eq1.25}
 \end{equation}
where $s_{\sL_i}$ is a covariant constant section of the restriction of
$(L^k, ka)$ to a $k$-cBS fibre $\sL_i$, and the sum is over all $k$-cBS
fibres of $\pi$.

 \begin{digr}{Mirror digression II} If a spLag cycle $\sL\subset X$ on a
CY manifold $X$ is an $n$-torus, its global deformation space
$\sM^{[\sL]}$ as a spLag cycle has dimension $n$ (see Mirror
digression~I). Following Strominger, Yau and Zaslow \cite{SYZ}, it is
reasonable to suppose that there exists a compactification
$\ov{\sM^{[\sL]}}$ such that the universal family defines a fibration
 \begin{equation}\pi\colon X\to \ov{\sM^{[\sL]}}
 \label{eq1.26}
 \end{equation}
with $T^n$ as generic fibre.

Then a polarization $L$ of $X$ gives the prequantization data
$(X,\om,L,a)$, and a finite set
 \begin{equation}
 \{\pi^{-1}(m_1),\dots,\pi^{-1}(m_N)\}
 \label{eq1.27}
 \end{equation}
of spLag BS cycles. Of course, they are all spLag cBS cycles. 

Now, if we are lucky, we can construct the dual fibration
 \begin{equation}
 \pi'\colon X'=\Pic(X/\ov{\sM^{[\sL]}})\to \ov{\sM^{[\sL]}},
 \label{eq1.28}
 \end{equation}
with $(\pi')^{-1}(m)=\Pic(\pi^{-1}(m))$ the space of flat Hermitian
connections on the trivial line bundle, and we get the finite set
(\ref{eq1.27}) as the intersection of two Lagrangian cycles in $S'$. But
these cycles admit an additional structure -- that of branes or
supercycles.
 \end{digr}

 \begin{digr}{Digression on branes} It is technically very convenient to
give any Lagrangian cycle $\sL$ the structure of a brane (or supercycle).
A {\em brane} or {\em supercycle} is a pair $(\sL,a)$, where $\sL$ is a
Lagrangian cycle and $a$ a flat connection on the trivial line bundle on
$\sL$. Then $X'$ is the {\em moduli space of supercycles (branes)} and the
projection $\pi'$ sends a supercycle to its Lagrangian support. The fibre
of $\pi'$ over $\sL$ is the moduli space of supercycles with support
$\sL$. We see below that the moduli space $X'$ admits a CY manifold
structure, and according to the SYZ conjecture, it is the {\em mirror
reflection} of $X$.
 \end{digr}

Let $L_0$ be the trivial line bundle on $X$ with the trivial connection
$a_0$. By definition the fibration $\pi'$ (\ref{eq1.28}) admits the ``zero
section''
 \begin{equation}
 s_0=\GFT(L_0)\subset X', \quad\text{where}\quad
 s_0=\{(\sL,a)\bigm|\text{$a$ is trivial}\}.
 \label{eq1.29}
 \end{equation}
We view it as an $n$-cycle in $X'$. On the other hand, if a holo\-morphic
line bundle $L$ on $X$ admits a Hermitian connection $a_L$ with curvature
form proportional to $\om$ then the restriction satisfies
 \begin{equation}
{a_L}\rest{\sL}\in (\pi')^{-1}(\sL)
 \label{eq1.30}
 \end{equation}
and we get a section
 \begin{equation}
s_L\subset X', \quad\text{where}\quad s_L=\{(\sL, {a_L}\rest{\sL}) 
 \label{eq1.31}
 \end{equation}
of the fibration $\pi'$.

Thus the expected number of corrected Bohr--Sommerfeld cycles equals the
{\em intersection number} of $n$-cycles
 \begin{equation}
N=[s_0]\cdot[s_L]
 \label{eq1.32}
 \end{equation}
(here $s_0=s_{\Oh_X}$, of course).

To prove that this number equals the geometric intersection number
(\ref{eq1.27}), we must investigate the geometry of cycles of this type.
This is our aim in this paper. But as an illustration we start with a very
simple example.

 \subsubsection*{CY$_1$: elliptic curves} Let $C$ be an elliptic curve
with a complex orientation, with a given holo\-morphic 1-form $\theta$ and
a metric $g$ that give the universal cover $U$ and the tangent bundle $TC$
the standard constant Hermitian structures (that is, the Euclidean metric,
symplectic form and complex structure $I$). The K\"ahler form $2\pi i\om$
gives a polarization of degree 1. We fix a smooth decomposition
 \begin{equation}
C=S^1_+\times S^1_-=\U(1)_+\times\U(1)_-,
 \label{eq1.33}
 \end{equation}
of $C$ as a commutative group, where the final equality fixes the zero
$o\in C$. Let $L$ be a holomorphic line bundle with a holomorphic
structure given by a Hermitian connection
$a$ with curvature form $F_a=\om$, and
 \[
L=\Oh_C(o).
 \] 
The decomposition (\ref{eq1.33}) induces a decomposition
 \begin{equation}
H^1(C,\Z)=\Z_+\times\Z_-,
 \label{eq1.34}
 \end{equation}
and we have two generators
 \begin{equation}
 [e]\in\Z_+, \quad[s]\in\Z_- \quad\text{with}\quad
 [e]^2=[s]^2=0 \quad \text{and}\quad [e]\cdot[s]=1.
 \label{eq1.35}
 \end{equation}

Any smooth ``irreducible'' cycle in $C$ is the image $\fie(S^1)$ of a
smooth embedding $\fie\colon S^1\to C$, that is, a smooth circle $\sL$ in
$C$ (with primitive cohomology class $[\sL]\in H^1(C,\Z)$).

Now, over every point $p\in C$ the {\em Lagrangian Grassmannian} $\LaGr
p=\LaGr(TS_p)=\U(1)$ (see (\ref{eq1.1}--\ref{eq1.2})) is the circle of tangent
directions, and the {\em oriented Lagrangian Grassmannization} (\ref{eq1.1})
of
$TS$ is
 \begin{equation}
\LaGr(C)=C\times\U(1).
 \label{eq1.36}
 \end{equation}
Thus the phase map (\ref{eq1.13}) sends a point $p\in\sL$ to the tangent
direction of $\sL$ at this point. We get the following result.

 \begin{prop}
 \begin{enumerate}
 \item A cycle $\sL$ is spLag iff it is smooth and its universal cover 
 \begin{equation}
 U(\sL)\subset U(C)=\C
 \label{eq1.37}
 \end{equation}
is a line {\em with rational slope}.

 \item any spLag cycle is defined by its cohomology class up to
translation.
 \end{enumerate}
 \end{prop}

 \begin{pf} Indeed, the tangent direction of a spLag cycle must be
constant, and the slope must be rational to get a compact cycle from the
universal line by the exponential map. \qed \end{pf}

Thus every primitive spLag cycle defines its slope
 \begin{equation}
m_I(\sL)\in\PP^1_{\Q},
 \label{eq1.38}
 \end{equation}
and its complete family of translates gives the SYZ fibration
 \begin{equation}
 \pi\colon C\to S^1=\sM^{[\sL]},
 \label{eq1.39}
 \end{equation}
with fibres spLag cycles in $[\sL]$, and base the moduli space of spLag
cycles in the fixed cohomology class.

Using the basis (\ref{eq1.35}), we can realize the basis $[e],[s]$ of the
cohomology lattice as families of spLag cycles with slopes
 \begin{equation}
m_I(s)=0 \quad \text{and} \quad m_I(e)=\infty.
 \label{eq1.40}
 \end{equation}
Then every primitive spLag cycle $\sL$ has slope
 \begin{equation}
 [\sL]=r[s]+d[e] \quad \text{with} \quad \hcf(r,d)=1 \quad
 \text{and} \quad m_I(\sL)=\frac{d}{r}\,.
 \label{eq1.41}
 \end{equation}
Finally, we remark that the intersection number of two spLag cycles is
determined by their slopes (that is, their cohomology classes):

 \begin{prop} Two different irreducible spLag cycles $[\sL_1]$ and
$[\sL_2]$ with slopes $d_i/r_i$ intersect transversally in 
 \begin{equation}
 \#\bigl(\sL_1 \cap\sL_2\bigr)=|d_1\cdot r_2-d_2\cdot r_1|
 \label{eq1.42}
 \end{equation}
points.

 \end{prop}

 \subsection*{Complex quantization} This is nothing other than {\em the
classical theory of theta functions}. Indeed, the decomposition
(\ref{eq1.34}) defines the collection of compatible {\em theta structures}
of every level
$k$: if $C_k\subset C$ is the subgroup of points of order $k$, the
decomposition (\ref{eq1.34}) defines a decomposition
 \begin{equation}
C_k=\Z_k^+\times\Z_k^-,
 \label{eq1.43}
 \end{equation}
and a decomposition of spaces of wave functions
 \begin{equation}
 \sH_{L^k}=H^0(C, L^k)=\bigoplus_{c\in\Z_k^-} \C\cdot\theta_c
 \quad\text{with}\quad
 \rk \sH_{L^k}=k,
 \label{eq1.44}
 \end{equation}
where $\theta_c$ is the theta function with {\em characteristic} $c$ (see
\cite{Mum}). 

 \subsection*{Real polarization} Realizing the class $[e]$ as a family
of spLag cycles gives us a real polarization
 \begin{equation}\pi\colon C\to S^1_-=\sM^{[e]}
 \label{eq1.45}
 \end{equation}
(see (\ref{eq1.39}) and (\ref{eq1.26})), and we can consider the dual fibration
 \begin{equation}\pi'\colon \Pic(C/S^1_-)\to S^1_- \quad
 \label{eq1.46}
 \end{equation}
with fibres
 \[
 (\pi')^{-1} (p)=\Hom(\pi_1(\pi^{-1}(p),\U(1)).
 \]
This fibration admits the section
 \begin{equation}
 s_0\in \Pic(C/S^1_-) \quad\text{with}\quad s_0 \cap (\pi')^{-1}
 (p)=\id\in\Hom(\pi_1(\pi^{-1}(p),\U(1)),
 \label{eq1.47}
 \end{equation}
so that we have a decomposition
 \begin{equation}
 \Pic(C/S^1_-)=(S^1)'\times S^1_-\,.
 \label{eq1.48}
 \end{equation}
Now we can lift the fibrations (\ref{eq1.45}) and (\ref{eq1.46}) to the
universal cover of
$S^1_-$:
 \begin{equation}
 \renewcommand{\arraystretch}{1.3}
 \begin{matrix}
 \C^*&\longrightarrow&C & \\
 \downarrow&& \vphantom{\scriptstyle{\pi}} \downarrow \scriptstyle{\pi}\\
 \R&\stackrel{\exp}{\longrightarrow} & S^1_-& \\
 \uparrow&& \vphantom{\scriptstyle{\pi'}}\uparrow \scriptstyle{\pi'}\\
 \C^*&\longrightarrow &C' & \\
 \end{matrix}
 \label{eq1.49}
 \end{equation}
where the vertical maps $\C^*\to\R$ send $z$ to $\log(|z|)$, and
 \begin{equation}
C=\C^* /\{q^n\} \quad \text{and} \quad C'=\C^* / \{|q|^2\cdot q^{-n}\}.
 \label{eq1.50}
 \end{equation}
Thus the 2-torus $C'$ is equipped with 
 \begin{enumerate}
 \item a symplectic form $\om'$;
 \item a complex structure $I'$.
 \end{enumerate}

 \paragraph{Important remark} We get an identification
 \begin{equation}
 H^1(C,\Z)=H^1(C',\Z)
 \label{eq1.51}
 \end{equation}
that extends to a smooth identification of $C$ and $C'$. But under this
identification, $C$ and $C'$ have {\em opposite orientations.} \medskip

Now we can apply geometric quantization to the real polarization
(\ref{eq1.45}) of the phase space $(C,\om, L^k,a_k)$, where $a_k$ is the
Hermitian connection defining the holomorphic structure on $L^k$. Sending
the line bundle $L^k$ to the character of the fundamental group of a fibre
gives the section
 \begin{equation}
s_{L^k}\subset C'=\Pic(C/S^1_-)
 \label{eq1.52}
 \end{equation}
and the Bohr--Sommerfeld subset of $S^1_-$ is 
 \[
s_0 \cap s_{L^k}\subset s_0=S^1_-\,.
 \]

 \begin{prop} The following hold on $C'$:
 \begin{enumerate}
 \item the cohomology classes of sections are
 \[
 [s_{L^k}]=[s_0]+k\cdot[e'],
 \]
where $[e']$ is the class of a fibre of $\pi'$;

 \item these sections, and the fibres of $\pi'$ are primitive spLag cycles
with slopes
 \[
 m_{I'}(s_0)=0, \quad m_{I'}(s_{L^k})=k
 \quad\text{and}\quad m_{I'}(e')=\infty;
 \]

 \item under the identification $s_0=S^1_-=\U(1)$, the intersection points
 \[
s_0 \cap s_{L^k}=\U(1)_k
 \]
are elements of order $k$ in $\U(1)$.
 \end{enumerate}
 \end{prop}

 \begin{pf} Consider the case $k=1$. Recall that after lifting to the
universal cover, a Hermitian connection $a$ satisfying $F_a=\om$ is
defined by a constant 1-form that gives a {\em linear} map of the
universal cover of $s$ to the universal cover of $e$. So we get (1), (2)
and (3) for $k=1$. For the general case we apply the isogeny
 \begin{equation}
 \la_k\colon s_0=\U(1)\to\U(1)\quad\text{given by}\quad \la_k(z)=z^k.
 \qed \label{eq1.53}
 \end{equation}
 \end{pf}

Therefore we get the decomposition
 \begin{equation}
 \sH_{\pi}^k=\bigoplus_{\rho\in\U(1)_k} \C\cdot s_{\rho}.
 \label{eq1.54}
 \end{equation}

 \begin{cor} The numerical quantization problem (see
(\ref{eq1.16}--\ref{eq1.18})) has an affirmative solution, that is,
 \[
 \rk \sH_{L^k}= \rk \sH_{\pi}^k.
 \]
 \end{cor}

Moreover, we get:

 \begin{prop} There exists an isomorphism
 \[
 \sH_{L^k}=\sH_{\pi}^k,
 \]
canonical up to a scaling factor.
 \end{prop}

 \paragraph{Proof--Construction:} We get this isomorphism up to the action
of the \hbox{$k$-torus} $(\C^*)^k$ (compare the decompositions
(\ref{eq1.44}) and (\ref{eq1.54})). But the canonical isomorphism is
defined by the action of the Heisenberg group $H_k$ on holomorphic
sections of the line bundle $L^k$ (by the theory of theta functions, see
\cite{Mum}) and the natural extension of the action of $H_k$ on the
collection of Bohr--Sommerfeld orbits. Both of these representations are
irreducible, thus the uniqueness of irreducible representations of $H_k$
gives a canonical identification of these spaces up to scaling.
 \qed\par\medskip

 \subsection*{Vector bundles of higher rank} The correspondence $L\to s_L$
can of course be extended to any holomorphic bundle $E\to s_E\subset C'$,
where $s_E$ is a multisection of the fibration $\pi'$ of degree $\rk E$.
This correspondence works perfectly in the case of arbitrary dimension. We
use it as frequently as the slant product in 4-dimensional gauge theory.
We call it the {\em geometric Fourier transformation}, or GFT. We use it in
any dimension (also for Fano varieties), but it is quite enough to
understand its action in the case of elliptic curves.

Every stable holomorphic vector bundle $E$ on $C$ carries a
Hermitian--Einstein connection $a_E$ that defines the holomorphic
structure on $E$ (we describe this in detail below). The curvature of
$a_E$ is
 \begin{equation}
 F_{a_E}=\deg E\cdot 2\pi i\cdot\om.
 \label{eq1.55}
 \end{equation}
Hence the restriction of $a_E$ to every fibre of $\pi$ is a flat Hermitian
connection on a circle, and thus
 \begin{equation}
 \begin{gathered}
 (a_E)\rest{\pi^{-1}(b)}=\chi_1\oplus\dots\oplus \chi_{\rk E},
 \quad\text{where}\\
 \chi_i\in\Hom(\pi_1(\pi^{-1}(b)),\U(1)))=\Pic(\pi^{-1}(b))=(\pi')^{-1}(b).
 \end{gathered}
 \label{eq1.56}
 \end{equation}
Thus we get the cycle
 \begin{equation}
 \GFT(E)\subset C',
 \label{eq1.57}
 \end{equation}
which is a {\em multisection} of the fibration $\pi'$ of degree equal to
$\rk E$. Of course, for a line bundle
 \begin{equation}
 \GFT(L^k)=s_{L^k}.
 \label{eq1.58}
 \end{equation}
Recall that $\pi'$ is a fibration in Abelian groups, so that we can add
all the points (\ref{eq1.56}) pointwise over $S^1_-$. We get a section
 \begin{equation}
 s_{\det E}=\GFT(\det E)\subset C'.
 \label{eq1.59}
 \end{equation} 

 \begin{prop}
 \begin{enumerate}
 \item $E$ is stable $\Longleftrightarrow$ $\GFT(E)$ is primitive.
 \item the multisection $\GFT(E)$ is a spLag cycle.
 \end{enumerate}
 \end{prop}

 \begin{pf} The statement (2) is local, and the proof for line bundles
works. To prove (1) we must use the inverse of GFT, which we describe
below. \qed \end{pf}

Now if two spLag cycles have {\em finite slopes} on $C'$ we can multiply
them fibrewise. Indeed, for any $b\in S^1_-$, any spLag cycle $s$ of
finite slope has finitely many intersections points $\{(s_i)b\}$ with this
fibre. We can define $\odot$-multiplication fibrewise:
 \begin{equation}
(s^1 \odot s^2)_b=\bigcup_{i,j} (s_i^1)_b\cdot (s_j^2)_b.
 \label{eq1.60}
 \end{equation}
 
 \begin{lem}
 \begin{enumerate}
 \item For any spLag cycles $s^1$ and $s^2$, the product $s^1 \odot s^2$
is the union of primitive spLag cycles;

 \item for any two holomorphic bundles $E_1$ and $E_2$ we have
 \[
 \GFT(E_1\tensor E_2)=\GFT(E_1) \odot \GFT(E_2).
 \]
 \end{enumerate}
 \end{lem}

The proof is immediate.

In \cite{A}, Atiyah constructed a collection of semistable indecomposable
bundles:
 \begin{equation}
\Oh_C=F_1, F_2,\dots,F_r=S^r F_2,\dots \quad\text{with $\rk F_i=i+1$,}
 \label{eq1.61}
 \end{equation}
uniquely determined by the condition $h^0(F_r)=1$. 

 \begin{prop}[\cite{A}] The free $\Z$-module $\sF$ generated by $[F_i]$ is
a commutative ring with multi\-pli\-cative structure given by tensor
product, and
 \begin{equation}
 \sF\iso\sR(sl_2(\C))
 \label{eq1.62}
 \end{equation}
-- the ring of finite dimensional representations of the algebra
$sl_2(\C)$.

 \end{prop}

For the collection (\ref{eq1.61}) of Atiyah bundles
 \[
 \GFT(F_r)=r\cdot s_0.
 \]

 \subsection*{The inverse problem} Every stable vector bundle $E$ on an
elliptic curve $C$ is uniquely determined up to translation by its
topological type $(r, d)$, where $r=\rk E$ and $d=\deg E$. Hence the
moduli space is
 \begin{equation}
 \sM_{r,d}^s=C.
 \label{eq1.63}
 \end{equation}
But every primitive spLag cycle $\sL$ is also uniquely determined up to
translation by its slope, but the moduli space
 \begin{equation}
 \sM^{[\sL]=(r, d)}=S^1_-\,.
 \label{eq1.64}
 \end{equation} 
Thus the map
 \begin{equation}
 \GFT\colon\sM_{r,d}^s\to\sM^{[\sL]=(r, d)}
 \label{eq1.65}
 \end{equation}
sending $E$ to $\GFT(E)$ is nothing other than the fibration
(\ref{eq1.39}). This map has fibres of dimension 1, and there is no chance
of recovering $E$ from $\GFT(E)$. However, this just means that we have
not used all of the information which gives $E$. Namely, we can define
$\GFT(E)$ as a supercycle (or brane) (see Brane digression). To see this,
consider the relative direct product
 \begin{equation}
C\times_{S^1_-} C'=\{(c, c')\in C\times C' \bigm|\pi(c)=\pi'(c')\},
 \label{eq1.66}
 \end{equation}
with its projections
 \begin{gather*}
 \pi'\colon C\times_{S^1_-} C'\to C, \quad
 \pi\colon C\times_{S^1_-} C'\to C' \\
 \text{and}\quad p\colon C\times_{S^1_-} C'\to S^1_-\,.
 \end{gather*}
Our cycle $\GFT(E))$ is embedded into $C\times_{S^1} C'$:
 \begin{equation}
 \GFT(E)=\{(c,c')\bigm|c\in s\text{ and }c'\in\GFT(E)\}\subset
C\times_{S^1_-}C'.
 \label{eq1.67}
 \end{equation}

Then the restriction of the pair 
 \begin{equation}
 (\pi')^*(E,a_E)\rest{\GFT(E)} 
 \label{eq1.68}
 \end{equation}
admits a tautological topologically trivial line subbundle $L$ with
Hermitian connection $s_{\tau}$, so that the pair 
 \begin{equation}
(\GFT(E),a_{\tau})=\sGFT(E)
 \label{eq1.69}
 \end{equation}
 is a supercycle. It is easy to check the following:

 \begin{prop}\label{pro4.7} Let\/ $\SM^{(r,d)}$ be the moduli space of
supercycles of slope\/ $d/r$. Then the map
 \[
 \sGFT\colon\sM^s_{r,d}\to \SM^{(r,d)}
 \]
is an isomorphism.
 \end{prop}

This result gives the solution of the inverse problem for GFT.

Finally, if an indecomposable holomorphic vector bundle $E$ is
topologically trivial, $E$ admits a {\em holomorphic flat connection}
$a_h$, and the restriction $(\pi')^*(E,a_E)\rest{\GFT(E)}$ defines a flat
connection on $\GFT(E)$, which may be non-Hermitian. Such a connection is
given by the conjugacy class $[M]$ of any matrix $M$. Thus
 \begin{equation}
 \sGFT(E)=(\GFT(E),[M]).
 \label{eq1.70}
 \end{equation}
Let $N_r$ be the standard $r\times r$ unipotent matrix. Then for the
collection of Atiyah vector bundles (\ref{eq1.61}), it is easy to see that
 \begin{equation}
 \sGFT(F_r)=(r\cdot s_0,[N_r])
 \label{eq1.71}
 \end{equation}
Thus the inverse problem has a positive solution for all semistable
holo\-morphic vector bundles on $C$.

 \subsection*{The K\"ahler-Hodge mirror map} To work out this
construction in the CY realm, we will need to bring to bear the whole
collection of tricks used to describe the ``Mirror symmetry'' phenomenon.
The construction for elliptic curves is only a model for the higher
dimensional case.

For an elliptic curve $C$, the algebraic cohomology ring is
 \begin{equation}
A_C=H^{2*}(C,\Z)=H^0(C,\Z)\oplus H^2(C,\Z).
 \label{eq1.72}
 \end{equation}
Alongside the multiplication in cohomology groups and the symmetric
bi\-linear form
 \begin{equation}
(u, v)=[u\cdot v]_2
 \label{eq1.73}
 \end{equation}
where $[w]_i$ denotes the $i$th homogeneous component of $w$, we have the
skew\-symmetric form
 \begin{equation}
\Span{u,v}=(u^*, v),
 \label{eq1.74}
 \end{equation}
where the operator 
 \begin{equation}
*\colon A_C\to A_C, \quad *(u_1,u_2)=(u_1,-u_2)
 \label{eq1.75}
 \end{equation}
sends the Chern character of a vector bundle to the Chern character of the
dual bundle. Thus we can view $A_C$ as the standard even unimodular
lattice of rank 2, and as the symplectic unimodular lattice of rank 2. 

Now consider the map
 \begin{equation}
 \begin{gathered}
 \mir\colon A_C\to H^1(C',\Z)\\
 \text{with}\quad \mir([C])=[s_0], \quad \mir([\pt])=[e'],
 \end{gathered}
 \label{eq1.76}
 \end{equation}
where $[C]$ is the fundamental class of $C$ and $[\pt]$ the class of a
point on $C$. Obviously
 \begin{equation} \mir(u)\cdot \mir(v)=\Span{u,v}=(u^*, v).
 \label{eq1.78}
 \end{equation}
Both of our curves are {\em marked}, that is, the bases $[s],[e]$ and
$[s_0],[e']$ are fixed. So the lattices $H^1(C,\Z)$ and $H^1(C',\Z)$ are
identified (see (\ref{eq1.51})). The moduli space of elliptic curves with
this type of marking is the classical {\em period domain}
 \begin{equation}
 \begin{aligned}
 \Om_{C'}&=\{(z_0,z_1) \bigm|\Im (z_1/z_0)>0\} \\
 & =\{\C\cdot\theta' \bigm|i\cdot\theta'\wedge\ov{\theta'}> 0\}\subset\PP
H^1(C', \C).
 \end{aligned}
 \label{eq1.79}
 \end{equation}

On the other hand, for any marked curve we have the {\em complexified
K\"ahler cone}:
 \begin{equation}
 \begin{gathered}
 \sK_C\subset H^2(C,\R)\oplus i\cdot H^2(C,\R)\\
 \text{given by}\quad \sK_C=\{u+i\cdot v \bigm|v>0\},
 \end{gathered}
 \label{eq1.80}
 \end{equation}
where the orientation on $H^2(C,\R)$ is given by the condition $\om>0$.

The complexification
 \begin{equation}
 \mir_{\C}\colon \sK_C\to \Om_{C'}
 \label{eq1.81}
 \end{equation}
of the map $\mir$ (\ref{eq1.76}) is called the {\em K\"ahler--Hodge mirror
map}. It is easy to see the following:

 \begin{prop} The K\"ahler--Hodge mirror map gives a holomorphic
isomorphism of the domains.
 \end{prop}

 \begin{cor} The K\"ahler--Hodge mirror map extends to an
identification:
 \begin{equation}
 \sK_C\times \Om_{C} \longleftrightarrow \sK_{C'}\times\Om_{C'}.
 \label{eq1.82}
 \end{equation}
 \end{cor}

Now we can consider the ring $A_C$ as the topological K-functor by sending
a vector bundle $E$ to its Chern character $(\rk E,\deg E)$. We get the map
 \begin{equation}
 \mir\colon K^0_{\topo} (C)\to H^1(C',\Z).
 \label{eq1.83}
 \end{equation}
Comparing this map with the map 
 \begin{equation}
 [\GFT]\colon K^0_{\topo} (C)\to H^1(C',\Z),
 \label{eq1.84}
 \end{equation}
sending a holomorphic bundle $E$ on $C$ to the cohomology class
$[\GFT(E)]$, we get the following:

 \begin{prop} For a holomorphic vector bundle $E$ on $C$
 \begin{equation}
 \mir(\ch(E))\cdot \sqrt{\td_C})=[\GFT(E)].
 \label{eq1.85}
 \end{equation}

 \end{prop}

Of course, here $\sqrt{\td_C}=(1, 0)$, so $u\cdot \sqrt{\td_C}=u$, but we
will see below that the shape of this formula is completely general in all
dimensions.

 \section{Application of geometric quantization to the CY$_2$ case}
 \markright{\qquad Geometric quantization and mirror symmetry}

Let $S$ be a K3 surface with complex orientation given by a holomorphic
\hbox{2-form} $\theta$, and a K\"ahler metric $g$ with K\"ahler form
$\om$. Recall that the target sphere $S^2_g$ of the complex phase map
(\ref{eq2.49}) parametrizes all almost complex structures compatible with
this metric. If our metric $g$ is Ricci flat, any almost complex structures
compatible with this metric is integrable. In this case, the complex
phase map is quite simple: for every 2-cycle $\Si$ of $S$, the complex
phase map $m\colon\Si\to S^2_g=\PP^1$ sends every point $p\in\Si$ to the
complex structure defined by $T\Si_p\subset TS$. Then $\Si$ is a sdAG
cycle iff its image is a point, corresponding to that complex structure on
$S$ for which our cycle $\Si$ is pseudoholomorphic.

Remark that the conformal class of our metric $g$ restricted to $\Si$
defines a complex structure $\Si_g$ on it. We can generalize the notion of
sdAG cycles: we say that a cycle $\Si$ is a {\em $g$-twistor} if its
complex phase map $m\colon \Si_g\to\PP^1$ is holomorphic. Of course the
constant map to a point is holomorphic.

{From} now on, we always start from an {\em algebraic} surfaces.

 \subsection*{The K\"ahler--Hodge mirror map} Consider the ring of
algebraic cycles
 \begin{equation}
 A_S\subset H^{2*}(S,\Z)=\bigoplus_{i=0}^2 H^{2i}(S,\Z).
 \label{eq2.1}
 \end{equation}
Every holomorphic vector bundle $E$ on $S$ defines the vector $\ch(E)$
which is its Chern character
 \begin{equation}
 \ch(E)=\bigl(\rk E,c_1(E),\half c_1(E)^2-c_2(E)\bigr).
 \label{eq2.2}
 \end{equation}
As in the CY$_1$ case (elliptic curves), alongside the standard symmetric
bilinear form
 \begin{equation}
 (u, v)=[u\cdot v]_4,
 \label{eq2.3}
 \end{equation}
induced by multiplication in cohomology, there is a second bilinear form,
which this time is again symmetric:
 \begin{equation}
\Span{u,v}=u^*\cdot v\cdot \td_S 
 \label{eq2.4}
 \end{equation}
where
 \[
 *\colon A_S\to A_S \quad\text{is given by} \quad
 (u_0,u_1,u_2)^*=(u_0,-u_1,u_2),
 \]
and $\td_S$ is the Todd class of our surface $S$:
 \begin{equation} 
\td_S=(1,0,2)\in A_S. 
 \label{eq2.5}
 \end{equation}
Moreover in our lattice $A_S$, the standard even unimodular sublattice of
rank~2
 \begin{equation}
 H=H^0(S,\Z)\oplus H^4(S,\Z)
 \label{eq2.6}
 \end{equation}
is distinguished such that
 \begin{equation}
 A_S=H\oplus \Pic S.
 \label{eq2.7}
 \end{equation}
Now the complete even dimensional cohomology ring $H^{2*}(X,\Z)$ contains
the sublattice of {\em transcendental} cycles
 \begin{equation} 
T_S=A_S^{\perp}\subset
H^{2*}(X,\Z).
 \label{eq2.8}
 \end{equation}

 \begin{mcdn}\label{mc2.1} We consider a K3 surface $S$ whose
transcendental lattice $T_S$ contains a hyperbolic lattice $H$; that is,
$T_S$ admits an orthogonal decomposition 
 \begin{equation}
 T_S=H\oplus \Pic S' \quad \text{where} \quad \Pic S'=H^{\perp}.
 \label{eq2.9}
 \end{equation}
(Here the final lattice is nothing but the orthogonal complement to $H$ in
$T_S$, but later, we realize it as the Picard lattice of the mirror
partner $S'$ of our surface $S$; this explains the notation.) Below, we
fix 
 \begin{enumerate}
 \item a decomposition (\ref{eq2.9});
 \item a basis $[e],[s]$ of $H$ such that
 \begin{equation}
 [e]^2=0, \quad [s]^2=-2, \quad[e]\cdot [s]=1.
 \label{eq2.10}
 \end{equation}
 \end{enumerate}
 \end{mcdn}

Moreover we suppose that the orthogonal decomposition over $\Q$
 \begin{equation}
 H^{2*}(S,\Z)=[H\oplus \Pic S]\oplus[H\oplus \Pic S'] 
 \label{eq2.11}
 \end{equation}
is defined over $\Z$.
Recall that each of these components $H$ has a basis: in the first,
 \begin{equation}
 [S]\in H^0(S,\Z), \quad[\pt]\in H^4(S,\Z),
 \label{eq2.12}
 \end{equation} 
with the multiplication table
 \begin{equation}
 [S]^2=0, \quad[\pt]^2=0, \quad[S]\cdot[\pt]=1, 
 \label{eq2.13}
 \end{equation}
and the second the basis (\ref{eq2.10}).

Of course we can change the basis (\ref{eq2.13}) to imitate an ``elliptic
pencil with a section'':
 \begin{equation}
 \begin{gathered} 
 {[\pt]}, \quad[S]-[\pt] \quad \text{so that}\\
 [\pt]\cdot ([S]-[\pt])=1, \quad ([S]-[\pt])^2=-2
 \end{gathered}
 \label{eq2.14}
 \end{equation}
just as for (\ref{eq2.10}) (this parallelism will be important below).

 \begin{digr}{Digression on algebraic geometry} In the algebraic geometric
situation, $A_S$ is called the {\em Mukai lattice}; every holomorphic
vector bundle $E$ defines its {\em Mukai vector}
 \begin{equation}
 m(E)=\ch(E)
\sqrt{\td_S}, 
 \label{eq2.15}
 \end{equation}
where
 \begin{equation}
 \sqrt{\td_S}=(1,0,1). 
 \label{eq2.16}
 \end{equation}
is the {\em positive root} of the Todd class (\ref{eq2.5}).
Thus for any pair of holomorphic bundles $E_1$ and $E_2$, the
Riemann--Roch theorem gives
 \begin{equation}
\chi (\Hom (E_1, E_2))=\ch(E_1^*)\ch(E_2) \td_S=m(E^*_1)\cdot m(E_2).
 \label{eq2.17}
 \end{equation}
In particular, the dimension of the moduli space $\sM_E$ of {\em simple}
bundles $E$ of topological type $\ch(E)$ is given by the formula
 \begin{equation}
 \dim_{\C}\sM_E=2-m(E)\cdot m(E^*). 
 \label{eq2.18}
 \end{equation}
Therefore to change the standard quadratic form (\ref{eq2.3}) to the more
exotic (\ref{eq2.4}), we must change the basis (\ref{eq2.12}) of the
sublattice
$H$ to the basis (\ref{eq2.14}) that corresponds to an elliptic pencil:
 \begin{equation}
 [S]-[\pt]=(1,0,-1), \quad [\pt]=(0,0,1);
 \label{eq2.19}
 \end{equation}
this is because multiplication by $\sqrt{\td_S}$ sends this to the basis
(\ref{eq2.12}).

It is natural to extend the Hodge structure of $S$ to the full lattice
(\ref{eq2.11}), thus getting Hodge structures on all components of the
moduli space of vector bundles over $S$. See \cite{T2} for a general
approach in this spirit.
 \end{digr}

The sublattice
 \[
A_S\oplus \Pic S\oplus H=(\Pic S')^{\perp}
 \]
in (\ref{eq2.11}) contains two sublattices $A_S$ and $\Pic S\oplus H$
which can be identified by the homomorphism
 \begin{equation}
 \begin{gathered}
 \mir\colon A_S\to\Pic S\oplus H\quad\text{defined by} \\
 \mir([S]-[\pt])=[s],\quad\mir([\pt])=[e], \\
 \text{and $\mir(L)=L$ for $L\in\Pic S$.}
 \end{gathered}
 \label{eq2.20}
 \end{equation}
Now consider the composite:
 \begin{equation}
 \begin{gathered}
 \ch\circ\mir\circ\,[i]\colon\Pic S\to H\oplus i\cdot\Pic S \\
 \ch\circ\mir\circ\,[i](L)=[s]+i\cdot c_1(L)-\half c_1(L)^2\cdot[e],
 \end{gathered}
 \label{eq2.21}
 \end{equation}
where
 \[
 [i]\colon A_S\to\bigoplus_{n=0}^2i\cdot H^{2n}(S,\Z) \quad\text{is
given by}\quad [i](u_0,u_1,u_2)=(u_0,i\cdot u_1,-u_2).
 \]

It is convenient to apply this map to the second part of the decomposition
(\ref{eq2.11}), that is, to $\Pic S'$ and to get the map
 \[
 \Pic S'\to H\oplus i\cdot \Pic S'
 \]

 \begin{prop} The complexification of this map gives the
K\"ahler--Hodge mirror map $\mir_{\KH}$ for K3 surfaces. 
 \end{prop}

It only remains for us to explain what it is. The K\"ahler--Hodge mirror
symmetry for K3 surfaces has a rich history; we refer to Dolgachev's survey
\cite{D} for the development of this subject. Recall that if we fix an
identification
 \begin{equation}
 H^{2*}(S,\Z)=\bigoplus_{i=0}^2 H^{2i}(S,\Z)=
 H\oplus H^{\oplus3}\oplus(-E_{8})^{\oplus2}
 \label{eq2.22}
 \end{equation}
(a surface with such identification is called a {\em marked\/} K3), we get
embeddings
 \begin{equation}
 A_S=A\subset H\oplus H^{\oplus3}\oplus(-E_{8})^{\oplus2},
 \label{eq2.23}
 \end{equation}
and
 \begin{equation} T_S=T\subset H^{\oplus3}\oplus(-E_{8})^{\oplus2}.
 \label{eq2.24}
 \end{equation}
Then in the projective space
 \begin{equation}
 \PP T\tensor \C 
 \label{eq2.25}
 \end{equation}
we get the ``intersection form quadric''
 \begin{equation}
 q= \bigl\{u\in T_S\tensor\C \bigm|u^2=0 \bigr\},
 \label{eq2.26}
 \end{equation}
and in it the {\em period domain}
 \begin{equation}
 \Om_A=\bigl\{u\in q \bigm|u\cdot \ubar>0 \bigr\}. 
 \label{eq2.27}
 \end{equation}
Following Nikulin, we call this the {\em period domain of marked K3s with
$A\subset\Pic$} (or $A$-condition on the Picard lattice). It parametrizes
marked surfaces $S$ with the condition
 \begin{equation}
 A/H\subset \Pic S. 
 \label{eq2.28}
 \end{equation}
On the other hand, the complexification of the sublattice
 \begin{equation}
 A/H\subset
 H^{\oplus3}\oplus(-E_{8})^{\oplus2}
 \label{eq2.29}
 \end{equation}
contains the {\em complexified K\"ahler cone}
 \begin{equation} 
 \sK_{A/H}=\bigl\{u+i\cdot v \bigm|v^2>0\bigr\}\subset \Pic
S\tensor\R\oplus i\cdot \Pic S\tensor\R. 
 \label{eq2.30}
 \end{equation}
Now assuming that our Main condition~1.1 holds (see
(\ref{eq2.9}--\ref{eq2.10})), we consider the complexified K\"ahler cone
 \begin{equation}
 \sK_{\Pic S'} 
 \label{eq2.31}
 \end{equation}
and the period domain $\Om_A$ (\ref{eq2.27}) of marked K3 surfaces with
$A\subset\Pic$. Then Proposition 2.1 says that the
complexification
$(\ch\circ\mir\circ\,[i])_{\C}$ of the map (\ref{eq2.21}) gives the map
$\mir_{\KH}$ as the composite 
 \begin{equation}
 \sK_{\Pic S'}
 \xrightarrow{(\ch \circ\mir\circ\,[i])_{\C}}
(T\tensor \C \setminus \{0\}) \stackrel{/\C^*}{\longrightarrow} \Om_A,
 \label{eq2.32}
 \end{equation}
and this map is {\em biholomorphic}. Note that the final $\C^*$-bundle
over $\Om_A$ is the principal bundle of the line bundle with fibres
$H^{2,0}$ for every K3 in this family.

Thus we get a topological identification of the two marked 4-manifolds $S$
and $S'$ underlying the two K3s. The two $H^2$ lattices of these are
identified (see (\ref{eq2.11})) and by Kulikov's theorem on surjectivity of
Torelli, the map $\mir_{\KH}$ extends to an identification
 \begin{equation}
 \sK_{\Pic S}\times \Om_{A_S} \equiv \sK_{\Pic S'}\times \Om_{A_{S'}},
 \label{eq2.33}
 \end{equation}
where $A_{S'}=H\oplus\Pic S'$.

Thus, for any integral class $[\om_I]$ such that $\om_I$ is the K\"ahler
form of some K\"ahler metric on $S_I$, the pair $([\om_I], S_I)$ acquires a
``mirror'' pair $([\om_{I'}], (S')_{I'})$. 

 But we get more from $\mir_{\KH}$: namely, consider the first part of the
composite (\ref{eq2.32}) without dividing out by complex scaling. Then we
get the triple 
 \begin{equation}
 \bigl([\om_{I'}], (S')_{I'},\theta'\bigr)
 \label{eq2.34}
 \end{equation}
where $\theta'$ is the complex orientation of $(S')_{I'}$. Recall that for
any triple of symplectic forms $(\Re\theta',\Im\theta',\om)$, we use
the usual normalization
 \begin{equation}
 \begin{gathered}
 (\Re\theta')^2=(\Im\theta')^2=\om^2; \\
 \Re\theta'\wedge\Im\theta'=\Re\theta'\wedge\om=
 \Im\theta'\wedge\om=0.
 \end{gathered}
 \label{eq2.35}
 \end{equation}
This complex orientation $\theta'$ is defined up to a rotation $e^{i\fie}$
and we have two closed symplectic forms 
 \begin{equation}
 \Re\theta' \quad \text{and} \quad \Im\theta'.
 \label{eq2.36}
 \end{equation}
Now using these relations, the cohomology class $[\om_{I'}]$ can be lifted
uniquely to the K\"ahler form $\om_{I'}$ of a K\"ahler metric $g'$ on
$S'$. Thus we get a K\"ahler structure on $S'$.

Moreover, we get an identification of 
 \begin{enumerate}
 \item the target sphere $S^2_g$ of $(i\cdot\om_I,S_I)$ and the target
sphere $S^2_{g'}$ of $(\om_{I'},(S')_{I'})$ and

 \item the circles $S^2\cap H\tensor\R$ and $(S^2)' \cap H\tensor\R$ in
such a way that the basis $[e],[s]$ goes to $[e'],[s_0]$.
 \end{enumerate}

 \subsection*{Remark on QFT} This K\"ahler-Hodge mirror symmetry can be
described in terms of the local structure of the moduli space of {\em
conformal fields theories} (see Aspinwall and Morrison \cite{AM}). In
particular, the orthogonal decomposition (\ref{eq2.22}) is nothing other
than the decomposition (\ref{eq2.11}) of \cite{AM}, and so on. 

 \subsection*{SYZ version of mirror symmetry for K3 surfaces}
Following the main idea of \cite{SYZ}, we want to realize all distinguished
cohomology classes by geometric objects such as spLag and sdAG cycles.

For this, we fix a {\em hyperK\"ahler metric} $g$ (that is, a Ricci
flat metric). The family of compatible complex structures and the
corresponding family of K\"ahler forms admits two distinguished algebraic
complex structures and K\"ahler (Hodge) form
 \begin{equation}
 I\in S^2_g ;
 \quad\om_I\in S^2_g 
 \label{eq2.37}
 \end{equation}
and the antipodal point $\Ibar$. This sphere $S^2_g$ is identified to the
target sphere of the complex phase map and we view it as the ordinary
sphere in $\R^3$:
 \begin{equation}
 S^2_g\subset\R^3=\Span{[\om_I],[\Re\theta],[\Im\theta]}\subset
H^2(S,\R),
 \label{eq2.38}
 \end{equation}
where $\theta$ is any complex orientation of $S_I$ with the normalization
(\ref{eq2.35}). Moreover,
 \begin{equation}
 [\om_I]\in \Pic S_I=A_{S_I} / H\in H^2(S,\Z)
 \label{eq2.39}
 \end{equation} 
(see (\ref{eq2.11})).

The three forms (\ref{eq2.38}) are pairwise orthogonal and define a basis
of
$\R^3$, viewed as the subspace of imaginary quaternions in the
4-dimensional space of quaternions $\HH$. Using this interpretation, we
relabel these forms
 \begin{equation}
 [\om_I],[\om_J],[\om_K], \quad \text{with $I,J,K\in\HH$}
 \label{eq2.40}
 \end{equation}
the standard basis of imaginary quaternions.
 \subsection*{Complex quantization} The algebraic class (\ref{eq2.39})
defines an
$I$-holomorphic line bundle $L$ on the complex surface $S_I$ with first
Chern class
 \begin{equation}
 c_1(L)=[\om_I]. 
 \label{eq2.41}
 \end{equation}
As any holomorphic bundle, $L$ carries a Hermitian structure and a unitary
connection $a_{L}$ given by the covariant derivative 
 \begin{equation}
 \nabla_L=\dd+\partial_I\fie
 \label{eq2.42}
 \end{equation}
 with curvature form 
 \begin{equation}
 F_{a_L}=2\pi i\cdot\om_I=\partial_I\dbar_I\fie,
 \label{eq2.43}
 \end{equation}
where $\fie$ is the potential function of our K\"ahler metric.

 \begin{rmk} Recall that on a simply connected 4-manifold $S$, any gauge
orbit is uniquely determined by its curvature form. \end{rmk}

Now we can consider the pair $(S,\om)$ as the phase space of a classical
mechanical system, and the full collection
 \begin{equation}
 (S,\om_I,L,a_L) 
 \label{eq2.44}
 \end{equation}
as its prequantization. Then the complex version of quantization provides
the space of wave functions
 \begin{equation}
 \sH_L=H^0(S_I, L).
 \label{eq2.45}
 \end{equation}
Our case is a priori different from the $\Pic=\Z$ case, and we need to
consider powers $L^k$ along with all the other integral purely imaginary
point of the K\"ahler cone (\ref{eq2.30}).

 \subsection*{Special real polarizations} We now construct a real
polarization of the phase space (\ref{eq2.44}), that is, a special
Lagrangian fibration
 \begin{equation}
 S\to S^2 
 \label{eq2.46}
 \end{equation}
in order to be able to apply our previous construction comparing the
results of quantizations in the complex and real polarizations of our
system (\ref{eq2.44}).

First of all, we have an orthogonal decomposition induced by (\ref{eq2.11})
 \begin{equation}
 H^2(S,\R)=
 (\Pic S_I\tensor\R) \oplus (H\tensor\R) \oplus (\Pic S'\tensor\R).
 \label{eq2.47}
 \end{equation}
The subspace $T_{S_I}\tensor\R$ contains the two planes
 \begin{equation}
 \Span{\Re\theta,\Im\theta} \quad \text{and} \quad
 H\tensor\R=\Span{[s],[e]}
 \label{eq2.48}
 \end{equation}
(see (\ref{eq2.9}--\ref{eq2.11})). It is easy to see that
 \[
 \Span{\Re\theta, \Im\theta}\ne \Span{[s],[e]}
 \]
(because the restrictions of the intersection form have different indexes).
On the other hand, $[s]$ and $[e]$ are orthogonal to $\om_I$, and thus
both have the same image under orthogonal projection to
$\R^3=\Span{\Re\theta,\Im\theta,\om}$.

Applying if necessary a {\em phase rotation}
 \begin{equation}
 \theta\to e^{i\cdot\fie}\cdot\theta,
 \label{eq2.49}
 \end{equation}
we can suppose that
 \begin{equation}
 J=\pr([e]) \quad \text{and} \quad K=\pr([s]).
 \label{eq2.50}
 \end{equation}
Then we have the following result:

 \begin{prop}
 \begin{enumerate}
 \item The cycle $[e]$ is algebraic in the complex structure $J$, and $[s]$
is algebraic in the complex structure $K$ (see (\ref{eq2.40})), that is
 \begin{equation}
 [e]\in\Pic S_J \quad \text{and} \quad [s]\in\Pic S_K.
 \label{eq2.51}
 \end{equation}

\item After a finite set of $-2$-reflections, and possibly a change of
basis, we can assume that\/ $[e]$ is the class of some irreducible elliptic
curve $e$ in the complex structure $J$, and $[s]$ the class of some
irreducible rational curve $s$ in the complex structure $K$.

\item The elliptic pencil\/ $|e|$ defines a $J$-holomorphic fibration
 \begin{equation}
 \pi: S_J\to\PP^1,
 \label{eq2.52}
 \end{equation}
that is, a spLag fibration over the compactification of the moduli space
 \begin{equation}
 \PP^1=\ov{\sM^{[e]}_J}
 \label{eq2.53}
 \end{equation}
of the complex phase $J$.

\item The irreducible smooth $K$-holomorphic curve $s$ is a section of this
pencil.
 \end{enumerate}
 \end{prop}

 \begin{pf} (2) follows from standard arguments Pyatetski\u\i-Shapiro and
Shafarevich, \cite[\S3]{PS}. The other statements follow directly from the
definitions and previous constructions. Remark that from the local
(twistor) description of complex directions at a point (see the previous
section), chosen complex directions intersect transversally and with the
same orientation.
 \qed \end{pf}

 \begin{rmk} We could consider these statements only for a {\em generic}
point of the period space $\Om_{A_S}$ (\ref{eq2.33}) to avoid various
technical problems. Then we can prove the same statements for the case
when our metric is not Ricci flat, but is an arbitrary K\"ahler metric. In
this case, our almost complex structures $J$ and $K$ are not integrable,
and our curves $e$ and $s$ are pseudoholomorphic. But the amazing fact is
that all the arguments of \cite[\S3]{PS} still work.
 \end{rmk}

We consider the fibration $\pi$ (\ref{eq2.52}) as a real polarization of
our phase space (\ref{eq2.44}), and investigate the Bohr--Sommerfeld
fibres of this system (see Definition~\ref{dfnBS}). We assume that the
degenerations in $\pi$ are as simple as possible, so that it has 24
singular fibres
 \begin{equation}
 \Sing(\pi)=\{p_1,\dots, p_{24}\}
 \label{eq2.54}
 \end{equation}
(below we
add arguments for the generic case). Let
 \begin{equation}
S_J^0=S_J\setminus\Sing(\pi) \quad \text{and} \quad\pi: S_J^0\to\PP^1 
 \label{eq2.55}
 \end{equation}
be a smooth fibration with fibres $T^2$ or $\C^*$. Consider the dual
fibration
 \begin{equation}
 \pi'\colon \Pic(S_J^0 /\PP^1)\to\PP^1. 
 \label{eq2.56}
 \end{equation}
This fibration has the ``zero section''
 \begin{equation}
 s_0\subset \Pic(S_J^0 /\PP^1)
 \label{eq2.57}
 \end{equation}
corresponding to trivial flat connections on the fibres. 

 \begin{prop}
 \begin{enumerate}
 \item This fibration in Abelian Lie groups can be compactified (by adding
24 points), to give a fibration
 \begin{equation}\pi'\colon S'\to\PP^1,
 \label{eq2.58}
 \end{equation}
where $S'$ is the underlying smooth 4-manifold of a K3 surface.

 \item The surface $S$ has the structure of principal homogeneous fibration
over the Jacobian fibration (\ref{eq2.56}).

 \item The section (\ref{eq2.57}) can be compactified to the smooth
oriented 2-cycle
 \[
 s_0\subset S'
 \]
which is topologically $S^2$. 
 \end{enumerate}
 \end{prop}

The lattice $H^2(S',\Z)$ contains a sublattice $H$ with fixed basis
 \begin{equation}
 \Span{[s_0],[e']}\subset H^2(S',\Z) 
 \label{eq2.59}
 \end{equation}
where $e'$ is a fibre of $\pi'$.

It is convenient to consider the orthogonal projection
 \begin{equation}
 \pr_H\colon H^2(S',\Z)\to H=\Span{[s_0],[e']}. 
 \label{eq2.60}
 \end{equation}
We are now ready to apply the geometric Fourier transformation for the aK
structure on $S'$ described above (see (\ref{eq2.33}--\ref{eq2.36})).

Moreover, we get an identification of the target sphere $S^2_g$ and the
target sphere $S^2_{g'}$ and, in the same vein, the circles $S^2\cap
H\tensor\R$ and $(S^2)' \cap H\tensor\R$, in such a way that the basis
$[e],[s]$ goes to $[e'],[s_0]$.

 \subsection*{GFT for holomorphic line bundles} We return to the set-up of
(\ref{eq2.41}--\ref{eq2.45}). The restriction of our pair $(L,a_L)$ to any
fibre $e_p$ of the fibration $\pi$ (\ref{eq2.52}) over $p\in\PP^1$ defines
a flat connection and a character
 \[
\bigl(\chi_p\colon\pi_1(e_p)\to\U(1)\bigr)\in\Pic e_p \quad
\text{and}\quad \chi_p\in e'_p\subset S'\setminus\Sing(\pi),
 \]
see (\ref{eq2.56}).

For any 1-cycle $\ga\in H_1(e_p,\Z)$, we have
 \begin{equation}
 \chi_p(\ga)=\exp\Bigl(2\pi i\int_{\ga}\partial_I\fie\Bigr)
 \label{eq2.61}
 \end{equation}
(see (\ref{eq2.42}--\ref{eq2.43})). Moreover, there exists a disc $D\subset
S$ such that
$\ga=\partial D$, and 
 \begin{equation}
 \chi_p (\ga)=\exp\Bigl(2\pi i\int_D \dbar_I \partial_I\fie\Bigr)=
 \exp\Bigl(2\pi i\int_{D}\om\Bigr)
 \label{eq2.62}
 \end{equation}
by Stokes' formula, and this number does not depend on the choice of $D$.
Indeed, the difference
 \begin{equation}
D-D'\in H^2(S,\Z) \implies\int_{D}\om-\int_{D'}\om\in\Z.
 \label{eq2.63}
 \end{equation}
Moreover, every 1-cycle $\ga\in H_1(e_p,\Z)$ defines a map
 \begin{equation}
\psi_{\ga}\colon T_S\tensor\R\to e'_p
 \label{eq2.64}
 \end{equation}
as follows: suppose that $D=D_I$ is a holomorphic disc in the complex
structure $I$. Then for any 2-form $q\in T_S\tensor\R$, the integral
 \begin{equation}
 \exp\Bigl(2\pi i\int_{D_I} q\Bigr)
 \label{eq2.65}
 \end{equation}
does not depend on the choice of $D_I$, because $D_I-(D_I)'$ is an
algebraic cycle and the integral of $q$ along it is 0.

Now we must consider a family of integrals of this form parametrized by
points $p\in\PP^1$, the base of our elliptic pencil $\pi$. Locally around
a regular point $p\in\PP^1$, we have special parameters of deformations of
spLag cycles.

 \subsection*{Special local coordinates on the moduli space of spLag
cycles} Let $u\subset\PP^1$ be a (very) small disc around a regular point
$p_0\in\PP^1$ and $U=\pi^{-1} (u)$. Then
 \begin{equation}\pi\colon U\to u=S^1_1\times S^1_2\times u
 \label{eq2.66}
 \end{equation}
is a smooth topologically trivial $T^2$-fibration over the disc, and we
have two fibred 3-manifolds
 \begin{equation}\pi\colon S^1_1\times u\to u \quad \text{and} \quad\pi\colon
S^1_2\times u\to u. 
 \label{eq2.67}
 \end{equation}
Now for every $p\in u$, the 1-cycles $(S^1_j)_p$ are contractible in $S$
for $j=1,2$, and for every such 1-cycle we can find an $I$-holomorphic disc
$D^j_p\subset S$ such that
 \begin{equation}
 \partial D^j_p=(S^1_j)_p.
 \label{eq2.68}
 \end{equation}
 Then
 \begin{equation}
 u_j(p)=\exp\Bigl(-2\pi\int_{D^j_p}\om\Bigr)
 \label{eq2.69}
 \end{equation}
are special coordinates on $u\subset\sM^{[e]}_J$ (see (\ref{eq2.52}),
Vafa \cite[formula (3.1)]{V} and the references given there). Moreover, if
for every $p\in u$
 \begin{equation}
 (S^1_{1})_p^*\times (S^1_{2})_p^*=\Hom(\pi_1(\pi^{-1}(p)),\U(1))\subset
(\pi')^{-1}(u)\subset S'
 \label{eq2.70}
 \end{equation}
is the dual torus with parameters 
 \begin{equation}
 \exp\bigl(-2\pi i\fie_j(p)\bigr),
 \label{eq2.71}
 \end{equation}
then
 \begin{equation}
 \wave{u}_j(p)=
 \exp\Bigl(-2\pi\Bigl(\int_{D^j_p}\om+i\fie_j(p)\Bigr)\Bigr)
 \label{eq2.72}
 \end{equation}
are local coordinates of the moduli space $\SM^{[e]}_J$ of spLag
supercycles \cite{V}. More precisely, we get special local coordinates of
 \begin{equation}
U'=(\pi')^{-1} (u)\subset S'=\SM^{[e]}_J.
 \label{eq2.73}
 \end{equation}

Returning to the character $\chi_p$ (\ref{eq2.62}), we get a map
 \begin{equation}
\psi_u\colon(\Pic S\oplus T_S\tensor\R)\times u\to U'
 \label{eq2.74}
 \end{equation}
given by the formula (\ref{eq2.65}) for any $q\in\Pic S\oplus
T_S\tensor\R$. This is well defined on $\Pic S$ because of (\ref{eq2.63})
and on $T_S\tensor\R$ because of the arguments after (\ref{eq2.65}).

In particular, for $q=\om$ we get a local section of the fibration $\pi'$:
 \begin{equation}
\psi_u (\om)=\GFT(L)_u\colon u\to U'.
 \label{eq2.75}
 \end{equation}
 
 \subsection*{Globalization} Now we can globalize these local complex
coordinates: it is easy to check the following result:

 \begin{prop}
 \begin{enumerate}
 \item For a form $\om$, we can glue the local complex coordinates
described above to give the complex structure on $S'\setminus\Sing\pi$
described by (\ref{eq2.33}--\ref{eq2.36});
 \item the local sections (\ref{eq2.75}) can be glued to a section
$\GFT(L)\subset S'$.
 \end{enumerate}
 \end{prop}

 \begin{warning} Although special coordinates give the same local
description of $\GFT(L)_u$ for any polarization $L$, the global complex
structure of $S'$ depends on the polarization $L$ (see (\ref{eq2.33})).
\end{warning}

 \begin{thm}
 \begin{enumerate}
 \item The image of\/ $[\GFT(L)]$ under the projection $\pr_H$
(\ref{eq2.60}) to the sublattice\/ $H$ is given by
 \begin{equation}
 \pr_H ([\GFT(L)])=[s_0]-\half c_1(L)^2[e'].
 \label{eq2.76}
 \end{equation}
 \item The cohomology class of\/ $[\GFT(L)]$ is
 \begin{equation} 
[\GFT(L)]=[s_0]+[\om_{I'}]-\half c_1(L)^2[e'] 
 \label{eq2.77}
 \end{equation}
(see (\ref{eq2.21})).

 \item $\GFT(L)$ is a spLag cycle of slope $-\half c_1(L) ^2$ with respect
to the aK structure $(\om_{I'}, (S')_{I'})$.
 \end{enumerate}
 \end{thm}

 \begin{rmk} The phase of a spLag cycle, that is, the point
 \begin{equation}
S^2_g \cap \Span{\Re\theta,\Im\theta}=\U(1)
 \label{eq2.78}
 \end{equation}
on the circle is determined as the slope or ratio of coordinates
$\pr([e])$ and $\pr([s])$ under the identification (\ref{eq2.50}), just as
in the elliptic curve case (see (\ref{eq1.41})).
\end{rmk}

 \begin{pf} The topological part of the proof can be deduced directly from
the multiplicative property of the Chern character, as in Cox and Zucker
\cite{CZ}, or in the more recent papers \cite{FMW}, \cite{BJPS},
transforming the additive structure of $H^2$ into the multiplicative
structure of sections of an elliptic pencil, together with (\ref{eq2.21}).
The proof that $\GFT(L)$ is spLag is purely local and is obvious in the
special coordinates (\ref{eq2.68}--\ref{eq2.71}). Finally, to compute the
slope we just need to use the map (\ref{eq2.74}). 
 \qed \end{pf}

 \begin{rmk} Instead of the local special coordinates
(\ref{eq2.66}--\ref{eq2.72}), we can use another trick: in fact
 \[
R^1\pi\Oh_S=T^*\PP^1,
 \]
admits the Eguchi--Hansen hyperK\"ahler metric $g_{EH}$ (the hyperK\"ahler
reduction of the standard $\U(1)$ action on $\C^2\times(\C^2)^*$) with the
sphere $S^2_{EH}$ of compatible complex structure, containing the two poles
$I,\Ibar$ of the ordinary and conjugate complex structures on $\PP^1$. For
any other complex structure,
 \[
p\colon (T^*\PP^1)_J\to\PP^1
 \]
is the affine bundle over the vector bundle $p\colon T^*\PP^1\to\PP^1$.
Now locally, the universal cover of $U$ (\ref{eq2.66}) is
 \[
 (\wave{\pi}\colon\wave{U}\to u)=p^{-1}(u),
 \]
and it carries the restriction of the Eguchi--Hansen metric to $p^{-1}(u)$,
the universal cover of $U$. The local lifting of $\GFT(L)$ is
pseudoholomorphic with slope in $S^2_{EH}$. \end{rmk}

Now it is reasonable to expect that there exists a finite set
 \begin{equation}
 \{m_1,\dots, m_{N_{L}}\}_{BS}\subset\PP^1
 \label{eq2.79}
 \end{equation}
of Bohr--Sommerfeld fibres with respect to $(L,a_L)$, so that we can
construct the space of wave functions for this real polarization $\pi$
(\ref{eq2.52}):
 \begin{equation}
 \sH_{\pi,L}=\sum_{\text{BS fibres}}\C\cdot s_{BS},
 \label{eq2.80}
 \end{equation}
that is, the direct sum of lines generated by covariant constant sections
on the Bohr--Sommerfeld fibres of $\pi$.

 \begin{thm}\label{th2.2} For every polarization $L$, we have equality
 \[
 \rk\sH_{L}=\rk\sH_{\pi,L}.
 \]
That is, the numerical quantization problem (see Remark following
(\ref{eq1.18})) has an affirmative solution.

 \end{thm}

 \begin{pf} By the Riemann--Roch theorem and the theorem on structure of
linear systems on a K3,
 \begin{equation}
 h^0 (L)=\half c_1(L)^2+2.
 \label{eq2.81}
 \end{equation}
On the other hand (\ref{eq2.79}) is the set of intersection points
 \begin{equation}
 \GFT(L) \cap s_0=\{m_1,\dots, m_{N_L}\}_{BS}. 
 \label{eq2.82}
 \end{equation}
Moreover, at any point of intersection of $\GFT(L)$ and $s_0$, the
orientations of the two sections correspond as for holomorphic curves in
the same complex structure. Thus the number of intersection points is
given by the formula
 \begin{equation}
 N_L=\bigl|[\GFT(L)]\cdot[s_0]\bigr|
 \label{eq2.83}
 \end{equation}
(compare (\ref{eq1.42}) for the case of elliptic curves). Finally, we should
remark that, as the zero section of a fibration in groups, the cycle $s_0$
has the opposite orientation. But then
 \begin{equation}
 N_L=\bigl((-[s_0]\bigr)\Bigl([s_0]-\half c_1(L)^2[e']\Bigr)=2+\half
c_1(L)^2,
 \label{eq2.84}
 \end{equation}
and we are done (see (\ref{eq2.81})). \qed \end{pf}

 \begin{rmk} It is natural to ask whether there is a ``natural
isomorphism'' between the spaces of wave functions of Theorem~\ref{th2.2};
that is, whether on an algebraic K3 surface whose mirror partner admits an
elliptic pencil, every complete linear system has some special ``theta
basis''.
 \end{rmk}

Now recall that the surface
 \[
 (S')_{I'}\setminus\Sing(\pi) \stackrel{\pi'}{\longrightarrow}\PP^1
 \]
is a fibration in groups, and for any polarization $L$ of $S_I$ we
obviously have
 \begin{equation}
 \GFT(L)\subset (S')_{I'}\setminus\Sing(\pi).
 \label{eq2.85}
 \end{equation}
Again, just as in the CY$_1$ case (see (\ref{eq1.60})) we can multiply
sections fibrewise
 \begin{equation}
 \GFT(L_1\tensor L_2)=\GFT(L_1) \odot \GFT(L_2).
 \label{eq2.86}
 \end{equation}
{From} this, we can extend GFT to all holomorphic line bundles on $S_I$.

Moreover, now in the same vein, for every holomorphic vector bundle $E$ on
$S_I$ that is stable with respect to our polarization $[\om]$, and such
that $c_1(E)=\la\cdot[\om]$, with $\la\in\Q$, the same construction (see
(\ref{eq1.55}--\ref{eq1.56})) with a Hermitian--Einstein connection defines a
cycle
 \begin{equation}
 \GFT(E)\subset S'
 \label{eq2.87}
 \end{equation}
which is a spLag multisection of $\pi'$ of degree $\rk E$, with class
given by a Chern character formula like (\ref{eq2.76}) and slope
determined from this formula as in (\ref{eq2.77}). We arrive once more to
the mirror reflection of the geometry of holomorphic bundles in terms of
the spLag and sdAG geometry under GFT.

 \subsection*{Inverse problem for the CY$_2$ case} Let $E$ be a vector
bundle on $S_I$ that is stable with respect to a polarization of $S_I$
given by an integral class $[\om]$, and with the condition
$c_1(E)=\la\cdot[\om]$, with $\la\in\Q$. Then $E$ carries a
Hermitian--Einstein connection with curvature form of the form
(\ref{eq1.55}). Therefore
 \begin{enumerate}
 \item this curvature defines the KH mirror pair $(S',\om')$.

 \item The topological type of $E$, given by the vector $m=m(E)$
(\ref{eq2.15}), together with $[\om']$, determines
 \begin{enumerate}
 \item a 2-dimensional cohomology class 
 \begin{equation}
 [\GFT(E)]\in H^2(S',\Z);
 \label{eq2.88}
 \end{equation} 

 \item a complex phase 
 \begin{equation}
 p=m_{I'} (\GFT(E))\in S^2_{g'},
 \label{eq2.89}
 \end{equation}
that is, an almost complex structure on $S'$ w.r.t.\ which the embedded
Riemann surface $\GFT(E)$ is pseudoholomorphic; and

 \item the complete linear system 
 \begin{equation}
 \bigl|[\GFT(E)]\bigr|^p=\sM^{[\GFT(E)]}
 \label{eq2.90}
 \end{equation}
(see \cite[(3.21)]{T3}) containing the given pseudoholomorphic curve
 \begin{equation}
 \GFT(E)\in\bigl|[\GFT(E)]\bigr|^p.
 \label{eq2.91}
 \end{equation}

 \end{enumerate}
 \item There is a map of the moduli space $\sM_m^s$ of stable vector
bundles over $S_I$ of topological type $m$:
 \begin{equation}
 \GFT\colon\sM_m^s\to\bigl|[\GFT(E)]\bigr|^p
 \label{eq2.92}
 \end{equation}
sending $E$ to $\GFT(E)$ (see (\ref{eq1.63}--\ref{eq1.65})).
 \end{enumerate}
As in the CY$_1$ case, we can give $\GFT(E)$ the structure of a supercycle
 \begin{equation}
 (\Si,\chi), \quad \text{for $\Si\in\bigl|[\GFT(E)]\bigr|^p$ and
 $\chi\in\Hom\bigl(\pi_1(\GFT(E)),\U(1)\bigr)$,}
 \label{eq2.93}
 \end{equation}
where the last space is the space of unitary flat connections on the
topologically trivial line bundle, that is, the space of classes of
characters of the fundamental group of $\Si$ into $\U(1)$.

Now recall that every pseudoholomorphic curve $\Si$ admits a complex
structure (see \cite[\S3]{T3}); thus we get an identification
 \begin{equation}
 \Hom\bigl(\pi_1(\GFT(E)),\U(1)\bigr)=J (\Si)
 \label{eq2.94}
 \end{equation}
with the Jacobian of $\Si$. Thus we can define the moduli space 
 \begin{equation}
 J(|\Si|^p)
 \label{eq2.95}
 \end{equation}
of such {\em pseudoholomorphic supercurves}.

Therefore our map GFT (\ref{eq2.92}) extends to a map
 \begin{equation}
 \sGFT\colon\sM_m^s\to J (|\Si|^p).
 \label{eq2.96}
 \end{equation}
Estimating the fibre of this map is the inverse GFT problem in the CY$_2$
case.

Now recall that the genus of any pseudoholomorphic cycles can be computed
by the adjunction formula \cite[\S3]{T3}, and the dimension of the moduli
space is given by (\ref{eq2.18}). {From} this and the formula for
$[\GFT(E)]$, we get
 \[
 \dim_{\R}\sM_m^s=4g(\GFT(E)).
 \]
Now using standard arguments of the theory of spectral curves, we can
prove the following result:

 \begin{prop} The map $\sGFT$ (\ref{eq2.96}) is of degree one. 
 \end{prop}

 \begin{rmk} It would be especially interesting to study the geometry of
the cycles $\sGFT(S^rT_S)$ in parallel with $\sGFT(S^rF_1)$ for the
collection of Atiyah bundles $S^rF_1$ over an elliptic curve discussed in
(\ref{eq1.61}); recall that in this case $S^rF_1$ is the bundle of $r$-jets.
 \end{rmk}

 We now consider the general geometric picture for the CY$_2$ case in terms
of the previous constructions:

 \begin{dfn} On any K\"ahler surface $(S,g)$, we define the {\em extended
Picard set\/} $\SDAG_g\subset H^2(S,\Z)$ to be the set of integral
combinations of classes of irreducible sdAG cycles of the same complex
phase:
 \begin{equation}
 \SDAG_g=\bigl\{[\Si]\in H^2(S,\Z) \bigm|\text{$\Si$ is sdAG}\bigr\}
 \label{eq2.97}
 \end{equation}
 \end{dfn}
Note that every cohomology class that can be realized as an irreducible
sdAG cycle admits two antipodal phase values.

 \begin{prop} If a metric $g$ on $S$ is hyperK\"ahler, that is, Ricci flat,
then any cohomology class $c\in H^2(S,\Z)$ admits only two antipodal
phases under any realization as a sdAG cycle.
 \end{prop}

 \begin{pf} The family of complex structures 
 \begin{equation}
 S_I\quad \text{for $I\in S^2$,}
 \label{eq2.98}
 \end{equation}
parametrized by the projectivization of the space of covariant constant
right chiral spinors is given by a conic in the period space. In this
space, the condition for a class to be algebraic is linear. Thus if some
hyperplane contains three (or more) points, it contains the whole conic. If
a cohomology class is algebraic in a complex structure $I$, it is also
algebraic in the complex structure $I'$. Thus if it is algebraic in one
complex structure, it is algebraic in every complex structure of the
family (\ref{eq2.78}), but this is forbidden by twistor theory, and we are
done. \qed \end{pf}

Now dividing the target sphere $S^2_g=\PP^1$ by the antipodal involution
gives us the real projective plane
 \begin{equation}
 \PP^2_{\R}=S^2/J,
 \label{eq2.99}
 \end{equation}
and a map
 \begin{equation}
 m_g\colon\SDAG\to\PP^2_{\R}.
 \label{eq2.100}
 \end{equation}
We get the subset
 \begin{equation}
 P_g=m_g(\SDAG_g)\subset\PP^2_{\R},
 \label{eq2.101}
 \end{equation}
of complex phase of sdAG cycles. Moreover, if 
 \[
 \sC=\bigl\{c\bigm|c^2>0\bigr\}\subset H^2(S,\R)
 \]
is the interior of the light cone, and the intersection
 \begin{equation}
 \SDAG\cap\sC\ne\emptyset,
 \label{eq2.102}
 \end{equation}
then this intersection is contained in one fibre of the map (\ref{eq2.76})
and defines an algebraic structure $(I,\Ibar)$, say the North and South
poles. Two such points on $S^2_g$ give special coordinates $(z_0,z_1)$ on
this complex Riemann sphere (don't forget the projective Hermitian
structure on the target sphere). Thus our projective plane admits a
structure of projective plane over $\Q$ with the set $\PP^2_{\Q}$ of
rational points as the set of rational directions. And the rational
projective line
$\PP^1_{\infty}$ parametrizes rational slopes on the equator. {From} this,
just as in the CY$_1$-case (see (\ref{eq1.34})), we get
 \begin{equation}
 P_g\subset S^2_{\Q},
 \label{eq2.103}
 \end{equation}
that is, every slope of a sdAG cycle is rational.

 \begin{rmk} This is the reason it was very convenient to start with the
case of an algebraic K3.
 \end{rmk}

The main condition~\ref{mc2.1} gives us two orthogonal points
 \begin{equation}
 m_g ([e]) \quad \text{and} \quad m_g ([s])
 \label{eq2.104}
 \end{equation}
 on the equator.

But the triple
 \begin{equation}
 (I, m_S([e]),m_S([s]))\subset P_S\subset S^2_{\Q}
 \label{eq2.105}
 \end{equation}
of pairwise orthogonal points switches on the $\GFT$ procedure in such a
way that at some time we get (\ref{eq2.80}) as an equality.

Now for our marked mirror pairs $(S,g)$ and $(S',g')$ (see
(\ref{eq2.33}--\ref{eq2.34})) with the identification
 \begin{equation}
H^2(S,\Z)=\Pic S_I\oplus H\oplus \Pic S'=H^2(S',\Z)
 \label{eq2.106}
 \end{equation}
(see (\ref{eq2.11})), we have
 \begin{equation}
 \SDAG_g\cap\SDAG_{g'}\subset H.
 \label{eq2.107}
 \end{equation}

 \section{GFT for the CY$_3$ case}
 \markright{\qquad Geometric quantization and mirror symmetry}

Recall that our main aim is to compare the geometry of holomorphic bundles
on algebraic (K\"ahler) varieties with the geometry of the images under GFT
of its (middle dimensional) cycles in the mirror partner. The mirror
partner is a priori a manifold of different topological type. But it is
natural in the CY case to expect that the mirror picture is symmetric, and
in particular, the mirror partner is again a CY manifold, and GFT can be
reversed. A conjectural construction of mirror symmetry admitting GFT was
proposed by Strominger, Yau and Zaslov in \cite{SYZ}, and developed in a
number of papers and preprints Gross and Wilson \cite{GW}, Gross
\cite{G1}, \cite{G2} and many others. A central technical problem of
realizing the SYZ program which is at present lacking is a description of
the global structure of the moduli space of spLag cycles. This is one of
the main differences between spLag geometry on the one hand and gauge
theory on the other. What is required is estimates as in Uhlenbeck's
theorem which would give an immediate description of compactifications of
the moduli spaces of connections and coherent sheaves as manifolds with
ends. Because of this Technical Mirror Asymmetry, the constructions
proposed below can be considered only as conditional, as in all the
current papers on the subject.

Thus, let $X$ be a simply connected CY 3-fold with a {\em complex
orientation}, that is, a holomorphic trivialization of $\det T^*X$. In the
CY$_3$ case, all the topological invariants of holomorphic vector bundles
and coherent sheaves on $X$ can again be realized as vectors in the ring
 \begin{equation}
 A_X=H^{2*}(X,\Q)=\bigoplus_{i=0}^n H^{2i}(X,\Q).
 \label{eq3.1}
 \end{equation}
This lattice carries the ordinary symmetric bilinear form, and the two
skew\-symmetric forms
 \begin{equation}
 (u,v)=[u\cdot v]_6=\sum_{i=0}^{3}u_i\cdot v_{3-i}
 \quad\text{and}\quad (u^*, v),
 \label{eq3.2}
 \end{equation}
where $u_i$ are the homogeneous components, and $*\colon A_X\to A_X$ is the
involution given by
 \[
 (u_0,\dots,u_3)^*=(u_0,-u_1,u_2,-u_3)
 \]
(see (\ref{eq1.73}--\ref{eq1.74})).

Alongside these forms there is {\em another skewsymmetric form}, given
by
 \begin{equation}
 \Span{u,v}=u^*\cdot v\cdot\td_X,
 \label{eq3.3}
 \end{equation}
where $\td_X$ is the Todd class of $X$: 
 \begin{equation}
 \td_X=\Bigl(1,0,\frac{1}{12}c_2(X),0\Bigr)\in A_X\quad\text{so that}\quad
(\td_X)^*=\td_X.
 \label{eq3.4}
 \end{equation}
We call
 \begin{equation}
c_2(X)=k_X\in H^4(X,\Z)
 \label{eq3.5}
 \end{equation}
the {\em second canonical class} of $X$.

The exotic bilinear form $\Span{u,v}$ (\ref{eq3.3}) can be reduced to the
standard bi\-linear form $(u^*, v)$ (\ref{eq3.2}) by multiplying by the
positive root
 \begin{equation}
 \sqrt{\td_X}=\Bigl(1,0,\frac{1}{24}k_X,0\Bigr)
 \label{eq3.6}
 \end{equation}
of the Todd class. Moreover, the topological type of any holomorphic
vector bundle $E$ is described by its {\em Mukai vector}
 \begin{equation}
 m(E)=\ch(E)\cdot \sqrt{\td_X}.
 \label{eq3.7}
 \end{equation}
Then for a pair of holomorphic vector bundles $E_1$ and $E_2$ on $X$,
formula (\ref{eq3.3}) is just the Riemann--Roch theorem:
 \begin{equation}
 \chi(E_1^*\tensor E_2)=\Span{\ch(E_1)^*,\ch(E_2)}=(m(E_1^*),m(E_2)),
 \label{eq3.7'}
 \end{equation}
and in particular the virtual dimension of the moduli space of an
$h$-stable vector bundle $E$ is equal to
 \begin{equation} 
 \vdim\sM_E= \Span{\ch(E),\ch(E)}=0
 \label{eq3.8}
 \end{equation}
because our exotic bilinear form (\ref{eq3.3}) is skewsymmetric. (For the
same reason, the virtual dimension of the moduli space of stable vector
bundles {\em with fixed determinant} on an elliptic curve is equal to 0
(see \S1).)

{From} this, we can define the integer CY $h$-invariant of the topological
class of an $h$-stable holomorphic vector bundle:
 \begin{equation}
 \CD_X^h (m)=\deg\sM^m_h=\# \{E\}
 \label{eq3.9}
 \end{equation}
to be the number of $h$-stable holomorphic vector bundles of topological
type $m=m(E)$ in the transversal case, and the top Chern class of the
obstruction bundle (under the deformation to the normal cone) in the
nontransversal case (see \cite{T1}). Recall the definition:

 \begin{dfn} The function
 \begin{equation}
 \CD_X^h \colon A_X\to\Z
 \label{eq3.10}
 \end{equation}
 is called the {\em Casson--Donaldson invariant} of $X$.
 \end{dfn}
This is the complex analog of the Casson invariant of a real 3-manifold,
and of the Donaldson polynomial of degree 0 (see Donaldson and Thomas
\cite{DT} and \cite{T1}). Remark that if some class in $a\in A_X$ cannot
be realized as the Mukai vector of a vector bundle, then $\CD_X^h (a)=0$.

Now the analog of the Picard lattice is the lattice
 \begin{equation}
 A_S / H=H^2(X,\Z) \oplus H^2(X,\Z)^*.
 \label{eq3.11}
 \end{equation}
of divisors and curves.

\subsection*{Complex polarization} For the line bundle $L=\Oh_X(h)$ of
any polarization $h$ we can consider our $X$ as the phase space of a
classical mechanical system with the complex pre\-quant\-iz\-a\-tion
structure
 \begin{equation}
 (X,\om,L,a_L),
 \label{eq3.12}
 \end{equation}
where $\om$ is the K\"ahler form, and $a_L$ the connection
(\ref{eq2.42}--\ref{eq2.43}) with curvature form $2\pi i\om$. Then we have
the space of wave functions
 \begin{equation}
 \sH_{\Oh_X(h)}=H^0(X, \Oh_X(h)) \quad \text{and} \quad X\subset \PP \sH^*.
 \label{eq3.13}
 \end{equation}

Now as in lower dimensions, there is a distinguished sublattice
 \begin{equation}
 H=H^0(X,\Z) \oplus H^6(X,\Z)\subset A_X
 \label{eq3.14}
 \end{equation}
with distinguished generators $[X]$ and $[\pt]$. The restriction of the
symmetric bilinear form $(\ ,\,)$ (\ref{eq3.2}) to this sublattice gives
the standard even unimodular lattice. But, whatever we might hope by
analogy with the case of K3 surfaces, this sublattice is not preserved by
multiplication by $\sqrt{\td_X}$ (\ref{eq3.7}), so that we obtain a new
standard sublattice
 \begin{equation}
H'=\sqrt{\td_X}\cdot H
 \label{eq3.15}
 \end{equation}
with a fixed basis $\sqrt{\td_X}$ and $[\pt]$, so that the form
(\ref{eq3.2}) restricted to it is again the standard even unimodular
lattice, but the restriction of the form $\Span{\ ,\,}$ (\ref{eq3.3})
gives the standard unimodular symplectic form.

\begin{rmk} A more invariant approach is to consider the rank~$3$
sublattice
 \[
 \Span{[X],k_X,[\pt]}.
 \]
Then both forms (\ref{eq3.2}) and (\ref{eq3.3}) restrict to degenerate
forms on this, with
$k_X$ generating the kernel of each. \end{rmk}

Now to switch on the mirror reflection of the holomorphic picture, we must
repeat step by step all the constructions proposed in the CY$_2$ case.
Recall that for any oriented supercycle $Y$ on $X$, we have the complex
phase map \cite[(6.16)]{T3}
 \[
 m_X\colon Y\to S^2
 \]
and $Y$ is a sdAG cycle if the image of $m_X$ is a point.

 \begin{dfn}
 \begin{enumerate}
 \item The subset 
 \[
 \SDAG=\bigl\{[Y]\in H^3(X,\Z) \bigm|\text{$Y$ is sdAG}\bigr\}\subset H^3(X,\Z)
 \]
is called the {\em Picard set} of a CY$_3$ manifold $X$.

 \item The symbol 
 \begin{equation}
P_X\subset S^2 / J=\PP^2_{\R}
 \label{eq3.16}
 \end{equation}
denotes the subset of complex phases of sdAG cycles.

 \item The symbol
 \begin{equation}
 F_X\subset\SDAG\times P_X
 \label{eq3.17}
 \end{equation}
denotes the subset of pairs $([Y], m_X (Y))$.
 \end{enumerate}
 \end{dfn}

The initial algebraic structure $(I,\Ibar)$, say, the North and South poles
of $S^2$, gives special coordinates $(z_0,z_1)$ on this complex Riemann
sphere, and our projective plane admits the structure of a projective
plane over $\Q$ with the set of rational points $\PP^2_{\Q}$ as the set of
rational directions; and the rational projective line $\PP^1_{\infty}$
parametrizes rational slopes on the equator.

Now following the SYZ program \cite{SYZ}, suppose that $P_X$ contains a
point $J=m_X(e)$ orthogonal to both poles, where $e=T^3$ is a 3-torus.
Then the SYZ conjecture predicts the existence of a fibration 
 \begin{equation}
 \pi\colon X\to
 \ov{\sM_{X}^{[e]}}=|e|, \quad \pi^{-1}(e_{\gen})=T^3,
 \label{eq3.18}
 \end{equation}
with special Lagrangian fibres over the compactification of the moduli
space -- the complete linear system of special Lagrangian cycles of the
cohomology class $[e]$ (see \cite{SYZ}).

\subsection*{Special real polarization} Suppose in addition that $P_X$
contains the point $K$ on the equator orthogonal to $m_X(e)$. Then there
exists a spLag cycle $s$ with slope $K=m_X(s)$ which is a section of the
fibration $\pi$ (\ref{eq3.18}). 

Thus the lattice $H^3(X,\Z)$ contains the distinguished sublattice
 \begin{equation}
H=\Z\cdot[s] \oplus\Z\cdot [e]\subset H^3(X,\Z).
 \label{eq3.19}
 \end{equation}

Passing to the relative Pic fibration (as in (\ref{eq2.56})) and using the
SYZ conjecture, we get the ``mirror'' fibration 
 \begin{equation}
\pi' \colon X'\to \ov{\sM_{X}^{[e]}}=|e|,
 \label{eq3.20}
 \end{equation}
and the complex structure on $X'$ is defined under the interpretation of
$X'$ as the moduli space $\ov{\SM_{X}^{[e]}}$ of special Lagrangian
supercycles. A more precise local description of this complex structure
is given below.

This fibration has a distinguished section 
 \begin{equation}
s_0\subset X'
 \label{eq3.21}
 \end{equation}
containing the supercycles with trivial flat connection. We again get a
distinguished skewsymmetric sublattice
 \begin{equation}
H=\Z\cdot [s_0] \oplus\Z\cdot [e']\subset H^3(X',\Z),
 \label{eq3.22}
 \end{equation}
where $[e']$ is the class of a fibre of $\pi'$ (\ref{eq3.20}).

\subsection*{Special local coordinates of the moduli space of splag
supercycles} Just as in the CY$_2$ case, let $u\subset
\ov{\sM_{X}^{[e]}}$ be a (very) small ball around a regular point
$p_0\in
\ov{\sM_{X}^{[e]}}$ and $U=\pi^{-1} (u)$. Then
 \[
\pi \colon U=S^1_1 \times S^1_2 \times S^1_3 \times u\to u 
 \]
is a smooth topologically trivial $T^3$-fibration over the ball and we
obtain three fibred 4-manifolds
 \begin{equation}
\pi \colon S^1_j \times u\to u \quad \text{where} \quad j=1, 2, 3 
 \label{eq3.23}
 \end{equation}
and three fibred 5-manifolds
 \begin{equation}
\pi \colon S^1_i \times S^1_j \times u\to u. 
 \label{eq3.24}
 \end{equation}

Now for every point $p\in u$, the 1-cycles $(S^1_j)_p$ are contractible in
$X$ for $j=1,2,3$, and for every such 1-cycle, we can find an
$I$-holomorphic disc $D^j_p\subset X$ such that
 \begin{equation}
\partial D^j_p=(S^1_j)_p.
 \label{eq3.25}
 \end{equation}
 Then
 \begin{equation}
 u_j(p)=\exp\Bigl(-2\pi\int_{D^j_p}\om\Bigr)
 \label{eq3.26}
 \end{equation}
are special coordinates on $u\subset\sM^{[e]}_J$.
Moreover, if for every $p\in u$,
 \begin{equation}
 (S^1_{1})_p^* \times (S^1_{2})_p^* \times (S^1_{3})_p^*=
 \Hom\bigl(\pi_1(\pi^{-1}(p)),\U(1)\bigr)\subset (\pi')^{-1}(u)\subset X'
 \label{eq3.27}
 \end{equation}
is the dual 3-torus with parameters $\exp\bigl(-2\pi i\fie_j(p)\bigr)$,
then
 \begin{equation}
 \wave{u}_j(p)=
 \exp\Bigl(-2\pi\Bigl(\int_{D^j_p}\om+i\fie_j(p)\Bigr)\Bigr)
 \label{eq3.28}
 \end{equation}
are local coordinates of the moduli space $\SM^{[e]}_J$ of spLag
super\-cycles. More precisely, we get special local coordinates on
 \begin{equation}
U'=(\pi')^{-1} (u)\subset S'=\SM^{[e]}_J.
 \label{eq3.29}
 \end{equation}
But in the same vein, we can get a system of dual coordinates: for every
point $p\in u$, the 2-cycles $(S^1_{i})\times(S^1_j)_p$ are contractible
in $X$ for $i\ne j=1,2,3$, and for every such 2-cycle we can find a
$I$-holomorphic surface $S^{i,j}_p\subset X$ such that
 \begin{equation}
 \partial S^{i,j}_p=
 \Bigl((S^1_{i})_p \times D^j_p\Bigr) \cup \Bigl(D^i_p \times 
 (S^1_j)_p\Bigr),
 \label{eq3.30}
 \end{equation}
where $D^j_p$ is as in (\ref{eq3.25}). Then
 \begin{equation}
u_{i,j} (p)=\exp\Bigl(-2\pi\int_{S^{i,j}_p} \om \wedge \om\Bigr)
 \label{eq3.31}
 \end{equation}
are {\em dual\/} special coordinates on $u\subset\sM^{[e]}_J$.

 \begin{rmk} Actually, we can take any disc $D^j_p$ and $4$-manifold
$S^{i,j}_p$ which are not $I$-holomorphic.
 \end{rmk}

Moreover, the 1-cycle $(S^1_j)_p\in H_1(e_p,\Z)$ defines a map
 \begin{equation}
\psi_j^1\colon H^2(X,\Z)\to e'_p
 \label{eq3.32}
 \end{equation}
by the formula 
 \[
\psi_j^1 (q)=\exp\Bigl(2\pi i\int_{D_p^j} q\Bigr),
 \]
and the 2-cycle $(S^1_{i})\times(S^1_j)_p\in H_2(e_p,\Z)$ defines a map
 \begin{equation}
 \psi_j^2\colon H^4(X,\Z)\to e'_p \quad\text{by}\quad
 \psi_j^2(q)=\exp\Bigl(\pi i\int_{S^{i,j}_p} q \Bigr).
 \label{eq3.33}
 \end{equation}
Composing these maps gives
 \begin{equation}
 \psi_u^j \colon H^{2j}(X,\Z) \times u\to U'
 \label{eq3.34}
 \end{equation}
given by these formulas for any $q\in H^{2j}(X,\Z)$. This map is well
defined by the same arguments as before ((\ref{eq2.63}) and
(\ref{eq2.65})).

In particular for $q=\om$, we get a local section of the fibration $\pi'$:
 \begin{equation}
 \psi_u^j (\om) \colon u\to U',
 \label{eq3.35}
 \end{equation}
and we can multiply these sections fibrewise (see (\ref{eq1.60})) to get a
local section
 \begin{equation}
 \psi_u (\om)=\psi_u^1(\om)\odot\psi_u^2(\om).
 \label{eq3.36}
 \end{equation}

\subsection*{Globalization} First of all, there exists a topological
globalization of the maps (\ref{eq3.34}) (see \cite{G1} and \cite{G2}).
That is, for any $q\in H^{2j}(X,\Z)$, the local sections $\psi_u^j(q)$
(\ref{eq3.35}) can be glued to a global section
 \begin{equation}
 \psi^j (q)\subset X'.
 \label{eq3.36'}
 \end{equation}
In particular we get the topological mirror map
 \begin{equation}
 \mir\colon A_X\to H^3(X',\Z), 
 \label{eq3.37}
 \end{equation}
given by
 \[
 \mir\Bigl((u_0,u_1,u_2,u_3)\sqrt{\td_X}^{-1}\Bigr)=
 u_0\cdot[s_0]+u_3\cdot[e']+[\psi^1(u_1)]+[\psi^2 (u_2)].
 \]
It is easy to check the following result.

 \begin{prop}
 \[
(\mir(u), \mir(v))= \Span{u,v},
 \]
where $(\ ,\,)$ is the intersection form in $H^3(X',\Z)$ and $\Span{\
,\,}$ the form (\ref{eq3.3}).
 \end{prop}

Hence for any $I$-polarization $[\om]$, and $L$ the holomorphic vector
bundle with $c_1(L)=[\om]$, we get the cycle
 \begin{equation}
\psi^j ([\om])=\GFT(L)\subset X'.
 \label{eq3.38}
 \end{equation}

 \begin{prop}
 \begin{enumerate}
 \item The intersection number $[\GFT(L)]\cdot[s_0]$ on $X'$ is given by
the formula
 \begin{equation}
([\GFT(L)])\cdot[s_0]=[\ch L\cdot\td_X]_6; 
 \label{eq3.39}
 \end{equation}
 \item the cycle $\GFT(L)$ is a spLag cycle of slope $\chi (L)$ equal to
the right-hand side of the Riemann--Roch equality.
 \end{enumerate}
 \end{prop}

The proof is again just the same as for the CY$_2$ case.

We can consider the fibration $\pi$ (\ref{eq3.18}) as a real polarization
of the system (\ref{eq3.12}), and geometric quantization provides the
space of wave functions
 \begin{equation}
\sH_{\pi}^L=\sum_{\text{$L$-BS fibres}} \C\cdot s_{\text{$L$-BS}},
 \label{eq3.40}
 \end{equation}
where the sum runs over all $L$-BS fibres of the system (\ref{eq3.40})
with respect to the real polarization $\pi$.

The number of such fibres is again the intersection number of
$s_0=\GFT(\Oh_X)$ and $\GFT(L)$. Using the same arguments as in the CY$_1$
and CY$_2$ cases, we get the following result.

 \begin{cor} For every polarization $L$, we have the equality
 \[
 \rk \sH_{L}=\rk\sH_{\pi, L} 
 \] 
(see (\ref{eq3.13}) and (\ref{eq3.40})). That is, the numerical
quantization problem (see Remark following (\ref{eq1.18})) has an
affirmative solution.
\end{cor}

Now in the same vein, for every holomorphic vector bundle $E$ on $X$ that
is stable with respect to our polarization $[\om]$ and such that
$c_1(E)=\la\cdot[\om]$ for $\la\in\Q$, the same construction with a
Hermitian--Einstein connection defines a cycle
 \begin{equation}
 \GFT(E)\subset X'
 \label{eq3.41}
 \end{equation}
which is a spLag multisection of $\pi'$ of degree $\rk E$, and its class is
given by the map $\mir$ (\ref{eq3.37}) and its slope is determined from
this formula as in (\ref{eq3.39}).

\subsection*{Inverse problem for the CY$_3$ case} Let $E$ be a stable
vector bundle on $X$ satisfying $c_1(E)=\la\cdot[\om]$ with
$\la\in \Q$. Then $E$ carries a Hermitian--Einstein connection, and its
curvature has the required shape. Hence:

 \begin{enumerate}
 \item This curvature defines the K\"ahler--Hodge mirror pair $(X,\om')$.

 \item The topological type of $E$, given by the vector $m=m(E)$
(\ref{eq3.7}), together with $[\om']$, determines a 3-dimensional
cohomology class 
 \begin{equation}
[\GFT(E)]\in H^3(X',\Z).
 \label{eq3.42}
 \end{equation} 

 \item It also determines the ``complete linear system'' 
 \begin{equation}
 \bigl|[\GFT(E)]\bigr|^p=\ov{\sM^{[\GFT(E)]}}
 \label{eq3.43}
 \end{equation} 
containing all spLag cycles of this cohomology class with the same slope
$p$.
 \end{enumerate}

Hence we again have a map of the moduli space $\sM_m^s$ of stable vector
bundles on $X$ of topological type $m$:
 \begin{equation}
\GFT\colon\sM_m^s\to \ov{\sM^{[\GFT(E)]}}
 \label{eq3.44}
 \end{equation}
sending $E$ to $\GFT(E)$.

Just as in the CY$_1$ and CY$_2$ cases, we can define a supercycle
structure on $\GFT(E)$ 
 \begin{equation}
 (\Si,\chi),\quad\text{where $\Si\in\bigl|[\GFT(E)]\bigr|^p$ and
 $\chi\in\Hom\bigl(\pi_1(\GFT(E)),\U(1)\bigr)$.}
 \label{eq3.45}
 \end{equation}
So we get the ``forgetful map''
 \begin{equation}
f \colon \ov{\SM^{[\GFT(E)]}}\to \ov{\sM^{[\GFT(E)]}}
 \label{eq3.46}
 \end{equation}
and the map
 \begin{equation}
 \sGFT\colon\sM_m^s\to \ov{\SM^{[\GFT(E)]}}
 \label{eq3.47}
 \end{equation}
in such a way that $\GFT=\sGFT\circ f$.

Estimating the fibre of this map is the inverse GFT problem in the CY$_3$
case. Using similar arguments as in the CY$_1$ and CY$_2$ cases, and by the
theory of deformations of spLag cycles, we get the following result.

 \begin{prop} The map $\sGFT$ (\ref{eq3.47}) is a degree $1$ map.
 \end{prop}

 \begin{pf} If $\sGFT(E_1)=\sGFT(E_2)$ then 
 \begin{equation}
 \sGFT(E_1^*\tensor E_2)=\sGFT(E_1^*)\odot\sGFT(E_2)
 \label{eq3.49}
 \end{equation}
contains the pair $([s_0],\chi=1)$ as a component. So there exists a
nontrivial homomorphism $E_1\to E_2$. {From} this by stability and standard
arguments we get $E_1=E_2$. \qed \end{pf}

Thus the situation is almost the same as in the CY$_1$ and CY$_2$ cases.

Now using the deformation theory of vector bundles on $X$ we get the
Kuranishi map
 \begin{equation}
K \colon H^1(X,\ad E)\to H^2(X,\ad E)
 \label{eq3.50}
 \end{equation}
such that locally, around the point $E\in\sM_m^s$, we have
 \begin{equation}
 \sM_m^s=K^{-1}(0).
 \label{eq3.51}
 \end{equation}

Thus we get the inequality
 \begin{equation}
 \dim_{\C} (K^{-1}(0))\le b_1(\GFT(E)).
 \label{eq3.52}
 \end{equation} 
because deformations of spLag cycles are unobstructed.

In particular, if deformations of $E$ are also unobstructed, then
 \begin{equation}
\rk H^1(X,\ad E)\le \rk H^1(\GFT(E), \C)
 \label{eq3.53}
 \end{equation}
which is very close to the equality sought by Vafa \cite[\S3]{V}. (Of
course, the differential of $\sGFT$ gives a homomorphism
 \begin{equation}
 \dd\sGFT\colon H^1(X,\ad E)\to H^1(\GFT(E), \C)
 \label{eq3.54}
 \end{equation}
which is an inclusion modulo the differential of the Kuranishi map.)

Furthermore, in many cases (and in the first instance, in the regular case
$\dim\sM_m^s=0$) it can be proved that (\ref{eq3.53}) is an equality.

The problem is to realize the generic spLag supercycle as $\sGFT(E)$
cycle. 

\subsection*{Example: $T_X$} It is well known that all deformations of the
rank 3 vector bundle $T_X$ are unobstructed (in the same way that
deformations of $X$ itself are unobstructed). Then using the SYZ
conjecture, one can prove that
 \begin{equation}
 H^1(X,\ad T_X)=H^1(\GFT(T_X), \C),
 \label{eq3.55}
 \end{equation}
the heart of the mirror conjecture.

 \section*{Acknowlegments}
 \markright{\qquad Geometric quantization and mirror symmetry}
 It will be clear to the reader that these constructions have been
developed and used by a number of physicists and mathematicians. I would
very much appreciate receiving comments and references from them; I would
also like to express my gratitude to Cumrun Vafa for his e-mail message.
Thanks are due to Miles Reid for tidying up the English.

\bigskip
\noindent
Andrei Tyurin, Algebra Section, Steklov Math Institute,\\
Ul.\ Gubkina 8, Moscow, GSP--1, 117966, Russia \\
e-mail: Tyurin@tyurin.mian.su {\em or} Tyurin@Maths.Warwick.Ac.UK\\
{\em or} Tyurin@mpim-bonn.mpg.de

\begin{thebibliography}{SYZ}
 \markright{\qquad Geometric quantization and mirror symmetry}

 \bibitem[AM]{AM} P. S. Aspinwall and D. R. Morrison, String theory on K3
surfaces, In: Mirror symmetry~II, A.M.S., 1997, pp.~703--716

 \bibitem[A]{A} M. F. Atiyah, Vector bundles over an elliptic curve, Proc.
London Math. Soc. (3) {\bf7} (1957), 415--452

 \bibitem[BJPS]{BJPS} M. Bershadsky, A. Johansen, T. Pantev, V. Sadov, On
four-dimensional compactifications of $F$-theory, Nucl. Phys. B{\bf505}
(1997), 165--201

 \bibitem[CZ]{CZ} D. A. Cox and S. Zucker, Intersection numbers of sections
of elliptic surfaces., Invent. math. {\bf53} (1979), 1--44

 \bibitem[D]{D} I. Dolgachev, Mirror symmetry for lattice polarized K3
surfaces, Alg. Geom. 4, J. Math. Sci {\bf81} (1996), 2599--2630

 \bibitem[DT]{DT} S. K. Donaldson and R. P. Thomas, Gauge theory in higher
dimensions, in The geometric universe (Oxford, 1996), OUP (1998),
pp.~31--47

 \bibitem[FMW]{FMW} R. Friedman, J. Morgan and E. Witten, Vector bundles
and $F$ theory., Comm. Math. Phys. {\bf187} (1997), 679--743

 \bibitem[G1]{G1} M. Gross, Special Lagrangian fibrations I: Topology, in
Integrable Systems and Algebraic Geometry, eds. Saito, Shimizu and Ueno,
World Scientific, 1998, pp. 156--193 (preprint alg-geom 9710006, 17 pp.)

 \bibitem[G2]{G2} M. Gross, Special Lagrangian fibrations II: to appear in
J. Diff. Geom., preprint alg-geom 9809073, 71 pp.

 \bibitem[GS]{GS} V. Guillemin and S. Sternberg, Symplectic Techniques in
Physics, CUP (1983)

 \bibitem[GW]{GW} M. Gross and P. M. H. Wilson, Mirror symmetry via 3-tori
for a class of Calabi-Yau threefolds, Math. Ann. {\bf309} (1997),
505--531

 \bibitem[H]{H} N. J. Hitchin, The moduli space of complex Lagrangian
submanifolds, preprint math/9901069, 18~pp.

 \bibitem[H2]{H2} N. J. Hitchin, The moduli space of special Lagrangian
submanifolds, Ann. Scuola Norm. Sup. Pisa Cl. Sci. (4) {\bf25} (1998),
503--515

 \bibitem[HL]{HL} R. Harvey, H. B. Lawson, Calibrated geometries, Acta
Math., {\bf148} (1982), 47--157

 \bibitem[Mum]{Mum} D. Mumford, Tata lectures on theta. I: Progr. Math,
28, Birkh\"auser (1983). II. Jacobian theta functions and differential
equations: Progr. Math, 43, Birkh\"auser (1984). III. Progr. Math, 97.
Birkh\"auser (1991)

 \bibitem[PS]{PS} I. I. Pjatecki\u\i -\v Sapiro and I. R. \v Safarevi\v c,
Torelli theorem for algebraic K3 surfaces, Izv. Akad. Nauk SSSR Ser. Mat.
{\bf35} (1971) 530--572 = Math. USSR Izvestiya {\bf5} (1971), 547--598

 \bibitem[S]{S} J. \'Sniatycki, Geometric quantization and quantum
mechanics, Applied Math Sciences, {\bf30} Springer (1980)

 \bibitem[SYZ]{SYZ} A. Strominger, S. T. Yau and E. Zaslow: Mirror symmetry
is T-duality, Nucl. Phys. B {\bf 479} (1996), 243--259

 \bibitem[T1]{T1} A. N. Tyurin, Non-Abelian analogues of Abel's theorem,
ICTP Trieste preprint (1997), 51~pp.

 \bibitem[T2]{T2} A. N. Tyurin, Symplectic structures on the moduli spaces
of vector bundles on algebraic surfaces with $p_g>0$, Izv. Akad. Nauk SSSR
Ser. Mat. {\bf52} (1988), 813--852 = Math. USSR-Izv. {\bf33} (1989),
139--177

 \bibitem[T3]{T3} Andrei Tyurin, Special Lagrangian geometry and slightly
deformed algebraic geometry (spLag and sdAG), Warwick preprint 22/1998,
alg-geom 9806006, 54~pp.

 \bibitem[V]{V} C. Vafa, Extending mirror conjecture to Calabi--Yau with
bundles, hep-th 9804131, 7~pp.

 \bibitem[W]{W} N. Woodhouse, Geometric Quantization, Oxford Math
Monographs, OUP (1980)

 \end{thebibliography}
\end{document}